\documentclass[11pt]{article}
\usepackage{amsfonts}
\usepackage{amssymb}
\usepackage{graphicx}
\usepackage{amsmath}
\usepackage{amsthm}
\usepackage{enumitem}
\usepackage{tikz}
\usepackage{mathrsfs}
\usepackage{cancel}
\usepackage{todonotes}
\usetikzlibrary{matrix}
\usetikzlibrary{arrows,shapes}
\usepackage{cite}
\usepackage{hyperref}
\hypersetup{colorlinks=true, urlcolor=blue}
\expandafter\let\expandafter\oldproof\csname\string\proof\endcsname
\let\oldendproof\endproof
\renewenvironment{proof}[1][\proofname]{
\oldproof[\ttfamily \scshape \bf #1. ]}{\oldendproof}
\newtheorem {theorem}{Theorem}[section]

\newtheorem {lemma}[theorem]{Lemma}
\newtheorem {proposition}[theorem]{Proposition}
\theoremstyle{definition}

\newtheorem{example}[theorem]{Example}
\newtheorem{remark}[theorem]{Remark}
\newtheorem{algorithm}[theorem]{Algorithm}

\numberwithin{equation}{section}
\normalsize
\makeatletter
\renewcommand\@biblabel[1]{#1.}
\makeatother
\def\Limsup_#1{\underset{#1}{\textrm {Lim sup }}}
\def\dom{{\rm dom\,}}
\def\epi{{\rm epi\,}}
\def\gph{{\rm gph\,}}

\def\ker{{\rm ker\,}}

\def\bd{{\rm bd\,}}

\def\inter{{\rm int\, }}
\def\dist{{\rm dist}}
\def\d{{\rm d}}
\def\argmin{{\rm argmin}}
\def\Q{{\mathcal{Q}}}
\def\nQ{{-\mathcal{Q}}}

\def\L{{\mathscr{L}}}

\def\B{{\mathbb{B}}}
\def\S{{\mathbb{S}}}
\def\R{\mathbb{R}}
\def\oR{\Bar{\R}}
\def\lm{\lambda}
\def\olm{\bar\lambda}

\def\ox{\bar{x}}
\def\oy{\bar{y}}

\def\ov{\bar{v}}

\def\ve{\varepsilon}
\def\Th{\Theta}
\def\Lm{{\Lambda}}
\def\tto{\rightrightarrows}
\def\ra{\rangle}
\def\la{\langle}
\def\ph{\varphi}
\def\dn{\downarrow}
\def\Bar{\overline}
\def\dd{\delta}
\def\sm{\hbox{${1\over 2}$}}
\def\O{\Omega}

\def\N{{\rm I\! N}}

\def\ri{\mbox{\rm ri}\,}
\def\x1k{{x^{k-1}}}
\def\xk{{x^k}}
\def\xkk{{x^{k+1}}}

\def\l1k{{\lambda^{k-1}}}
\def\lk{{\lambda^k}}
\def\lkk{{\lambda^{k+1}}}
\def\hlk{{\hat\lambda^k}}
\def\sig1k{{\sigma_{k-1}}}
\def\sigk{{\sigma_k}}
\def\sigkk{{\sigma_{k+1}}}
\def\ro1k{{\rho_{k-1}}}
\def\rok{{\rho_k}}

\def\e1k{{\varepsilon_{k-1}}}
\def\ek{{\varepsilon_k}}

\def\s1k{{p^{k-1}}}

\def\skk{{p^{k+1}}}

\def\alk{{\alpha^k}}
\def\C1k{{\mathcal{C}_{k-1}}}

\def\V1k{{V_{k-1}}}

\def\bx{{\bar x}}

\def\bl{{\bar\lambda}}
\def\hl{{\hat\lambda}}
\def\della{{\gamma_{\bar\lambda}}}

\def\ela{{\epsilon_{\bar\lambda}}}
\def\ellla{{\ell_{\bar\lambda}}}
\def\rhola{{\rho_{\bar\lambda}}}
\def\emp{\emptyset}
\def\hat{\widehat}
\def\tilde{\widetilde}
\def\epsilon{\varepsilon}

\begin{document}
\vspace*{0.5in}
\begin{center}
{\large\bf AUGMENTED LAGRANGIAN METHOD FOR SECOND-ORDER CONE PROGRAMS UNDER SECOND-ORDER SUFFICIENCY}\\[1ex]
\medskip

NGUYEN T. V. HANG\footnote{Department of Mathematics, Wayne State University, Detroit, MI 48202, USA and Institute of Mathematics, Vietnam Academy of Science and Technology, Hanoi 10307, Vietnam (hangnguyen@wayne.edu). Research of this author was partly supported by the USA National Science Foundation under grants DMS-1512846, DMS-1808978 and by the USA Air Force Office of Scientific Research under grant \#15RT04.},
BORIS S. MORDUKHOVICH\footnote{Department of Mathematics, Wayne State University, Detroit, MI 48202, USA (boris@math.wayne.edu). Research of this author was partly supported by the USA National Science Foundation under grants DMS-1512846 and DMS-1808978, by the USA Air Force Office of Scientific Research under grant \#15RT04, and by the Australian Research Council under Discovery Project DP-190100555.}
and M. EBRAHIM SARABI\footnote{Department of Mathematics, Miami University, Oxford, OH 45065, USA (sarabim@miamioh.edu).}
\end{center}

\noindent{\bf Abstract.} This paper addresses problems of second-order cone programming important in optimization theory and applications. The main attention is paid to the augmented Lagrangian method (ALM) for such problems considered in both exact and inexact forms. Using generalized differential tools of second-order variational analysis, we formulate the corresponding version of second-order sufficiency and use it to establish, among other results, the uniform second-order growth condition for the augmented Lagrangian. The latter allows us to justify the solvability of subproblems in the ALM and to prove the linear primal-dual convergence of this method.


\vspace*{-0.15in}

\section{Introduction}\label{intro}\vspace*{-0.05in}

In this paper we consider the class of constrained optimization problems belonging to {\em second-order cone programming} (SOCP) that are given in the form:
\begin{equation}\label{CP}
\mbox{minimize }\;f(x)\;\textrm{ subject to }\;\Phi(x)\in\Q,
\end{equation}
where both mappings $f\colon\mathbb{R}^n\to\mathbb{R}$ and $\Phi\colon\mathbb{R}^n\to\mathbb{R}^{m+1}$ are twice continuously differentiable ($\mathcal{C}^{2}$-smooth) around the reference points, and where the underlying set $\Q$ is the second-order (Lorentz, ice-cream) cone in $\mathbb{R}^{m+1}$ defined by
\begin{equation}\label{icecream}
\Q:=\big\{y=(y_0,y_r)\in\mathbb{R}\times\mathbb{R}^{m}\,\big|\,\Vert y_r\Vert\le y_0\big\}.
\end{equation}
Problems of this type (SOCPs for brevity) constitute a remarkable subclass of nonpolyhedral conic programs that has been well recognized in constrained optimization and various applications; see, e.g., \cite{ag,br,bs,mos} and the references therein.

Our main attention is paid to developing a numerical method to solve \eqref{CP} that involves the {\em augmented Lagrangian} $\L\colon\mathbb{R}^n\times\mathbb{R}^{m+1}\times(0,\infty)\to\mathbb{R}$ associated with this problem, which is defined  by
\begin{equation}\label{aL}
\L(x,\lambda,\rho):=f(x)+\dfrac{\rho}{2}\dist^2\big(\Phi(x)+\rho^{-1}\lambda;\Q\big)-\dfrac{1}{2}\rho^{-1}\Vert\lambda\Vert^2,\;\;(x,\lm,\rho)
\in\mathbb{R}^n\times\mathbb{R}^{m+1}\times(0,\infty),
\end{equation}
where $\lambda\in\R^{m+1}$ is a (vector) multiplier, and where $\rho>0$ is a penalty parameter of $\L$. The principal idea of the {\em augmented Lagrangian method} (ALM) for \eqref{CP} is to solve a sequence of unconstrained problems which objectives are defined by the augmented Lagrangian \eqref{aL} at a given multiplier-parameter pair $(\lm,\rho)$; namely,
\begin{equation}\label{aup}
\mbox{minimize }\;\L(x,\lm,\rho)\;\mbox{ over }\;x\in\R^n.
\end{equation}
This means that, given a multiplier $\lambda$ and a penalty parameter $\rho$, the ALM solves the unconstrained problem \eqref{aup} for the primal variable $x$ and uses the obtained value to update both the multiplier and penalty parameter in the next iteration.

The ALM was first proposed independently by Hestenes and Powell for nonlinear programming problems (NLPs) with equality constraints  \cite{Hes69,Pow69} and was originally known as the method of multipliers. For the latter framework, Powell observed in \cite{Pow69} that the ALM converges locally with an arbitrarily linear rate if one started the method with a sufficiently high penalty factor (but without the requirement of driving the penalty parameter to infinity) and from a point sufficiently close to a primal-dual pair that satisfies the standard second-order sufficient conditions (SOSC). This is an appealing feature of the ALM, since it provides a numerical stability that cannot be achieved in the usual smooth penalty method.

The ALM was largely extended to various settings of NLPs as well as convex programming with both equality and inequality constraints by Rockafellar \cite{Roc73a,Roc73b,Roc74}; see also the
monographs \cite{Ber82,NoW06,Rus06} and the references therein. Further theoretical and practical developments for this method were achieved in \cite{abm1,abm3,abm2,ah}  for NLPs. Furthermore, a
modified version of the ALM was suggested in \cite{abm1} to alleviate the possibility of unboundedness of dual sequences in the standard ALM; see \cite{bm14} for a detailed discussion on this
modification of the ALM and its successful implementation.

The classical results for the linear convergence of the ALM in NLP framework impose the SOSC, the linear independence constraint qualification (LICQ), and the strict complementarity condition,
which all together guarantee the uniqueness of the primal solution as well as the corresponding dual solution/multiplier.

More recently, the study of the ALM has been growing with important theoretical developments. On one hand, various attempts have been made to relax the restrictive assumptions for the convergence of this method in the NLP settings. In such a framework, Fern\'andez and Solodov achieved in \cite{FeS12} a remarkable progress for NLPs by proving that the linear convergence of the primal-dual sequence                                                                                                                              in the ALM can be ensured if the  SOSC alone is satisfied. This result significantly improved the classical ones for NLPs by verifying that neither the LICQ nor the strict complementarity condition is required for local convergence analysis of the ALM. A further improvement was obtained in Izmailov et al. \cite{IKS15} by showing that the conventional SOSC utilized in \cite{FeS12} can be replaced by the noncriticality of Lagrange multipliers for problems with equality constraints. On the other hand, the ALM has been studied for other major classes of constrained optimization including SOCPs \cite{lz} and semidefinite programming problems (SDPs) \cite{ssz}. For ${\cal C}^2$-cone reducible problems of conic programming (in the sense of Bonnans and Shapiro \cite{bs}), Kanzow and Steck \cite{KaS17,KaS19} established the linear convergence of the primal-dual sequence generated by modified versions of the ALM under the SOSC and strong Robinson constraint qualification; the latter yields that the Lagrange multiplier is unique. However, the solvability of subproblems in the ALM was not addressed in these papers. We also refer the reader to the paper by Cui et al. \cite{cst} and the bibliography therein for recent developments on the ALM for particular classes of convex composite problems of conic programming.\vspace*{0.03in}

The major goal of this paper is to develop both {\em exact} and {\em inexact} versions of the ALM for SOCPs under fairly {\em mild assumptions}. We aim first at establishing the {\em solvability} and {\em Lipschitzian stability} of the ALM subproblems by imposing merely the corresponding {\em SOSC} for \eqref{CP} in the general case of {\em nonunique} Lagrange multipliers. Having this, we verify a local {\em primal-dual convergence} of iterates with an arbitrary {\em linear rate} by assuming in addition the {\em uniqueness} of multipliers. Similarly to Fern\'andez and Solodov \cite{FeS12}, our approach revolves around the {\em second-order growth condition}   for the augmented Lagrangian \eqref{aL}. To the best of our knowledge, the origin of such a second-order growth condition for NLPs goes back to Rockafellar in \cite[Theorem~7.4]{r93} from which \cite{FeS12} significantly benefits. However, in contrast to \cite{r93}, \cite {FeS12} as well as to the vast majority of other publications on numerical optimization, we achieve our goal for \eqref{CP} by employing the concepts of the {\em second subderivative} and {\em twice epi-differentiability} of extended-real-valued functions in the framework of second-order variational analysis.\vspace*{0.03in}

The outline of the paper is as follows. In Section~\ref{sect02} we recall the notions of variational analysis and generalized differentiation with some {\em preliminary material} used in the paper.

The main result of Section~\ref{sect03}, which is of its own interest, provides an {\em error bound} estimate for the canonically perturbed KKT system associated with \eqref{CP} under the SOSC and a certain {\em calmness} property of the multiplier mapping with respect to perturbations that automatically holds for NLPs. We also present here an example showing that the imposed calmness property is essential for the validity of the error bound in the SOCP setting and then discuss efficient conditions ensuring the fulfillment of this calmness for nonpolyhedral SOCPs.

Section~\ref{paL} conducts a comprehensive {\em second-order variational analysis} of the augmented Lagrangian \eqref{aL} associated with the second-order cone program \eqref{CP}. Based on the obtained precise computation of the second subderivative of \eqref{aL}, we characterize here the {\em second-order growth} condition for \eqref{aL} via the SOSC and then establish its {\em uniform} counterpart needed in the general case of nonunique Lagrange multipliers.

The concluding Section~\ref{sect05} provides a detailed {\em solvability}, {\em stability}, and {\em local convergence analysis} of the suggested ALM algorithm for SOCPs that strongly exploits the SOSC and obtained second-order growth conditions. Our analysis includes the proof of solvability of the ALM subproblems in both exact and inexact versions and then establishes the linear convergence of primal-dual iterates to the designated solution of the KKT systems under the SOSC by using the established {\em robust isolated calmness} and {\em upper Lipschitzian} properties of the corresponding perturbed multiplier mappings. In this way we obtain explicit relationships between the constants involved in the algorithm and the imposed assumptions on the given data.\vspace*{0.03in}

Throughout the paper we use the standard notation and terminology of variational analysis and conic programming; see, e.g., \cite{bs,m18,RoW98}. Recall that $\B$ and $\S$ stand for the closed unit ball and unit sphere, respectively, of the space in question, that $\B_\gamma(x):=x+\gamma\B$  is the closed ball centered at $x$ with radius $\gamma>0$, that $\N:=\{1,2,\ldots\}$, and that $A^*$ indicates the matrix transposition. Given a nonempty set $\Omega\subset\R^n$, the symbols $\inter\Omega$, $\ri\O$, $\bd\Omega$, and $\Omega^\perp$ signify its interior, relative interior, boundary, and orthogonal complement space, respectively. The indicator function of $\O$ is defined by $\dd_\O(x):=0$ for $x\in\O$ and $\dd_\O(x):=\infty$ otherwise, dist$(x;\O)$ signifies the distance between $x\in\R^n$ and the set $\O$, and the projection of $x$ onto $\O$ is denoted by $\Pi_\O(x)$. As in \eqref{icecream}, we often decompose a vector $y\in\Q\subset\mathbb{R}^{m+1}$ into $y=(y_0,y_r)$ with $y_0\in\mathbb{R}$ and $y_r\in\mathbb{R}^{m}$. Taking this decomposition into account, denote $\widetilde y:=(-y_0,y_r)$. Similarly,   for the mapping  $\Phi\colon\R^n\to\R^{m+1}$ with $\Phi=(\Phi_0,\ldots,\Phi_{m})$, we  denote by $\widetilde\Phi(x)$ the vector $(-\Phi_0(x),\Phi_r(x))$ for any $x\in\R^n$.\vspace*{-0.15in}

\section{Preliminaries from Variational Analysis}\label{sect02}\vspace*{-0.05in}

We start this section with recalling those constructions of variational analysis and generalized differentiation, which are broadly employed in what follows; see \cite{bs,m18,RoW98} for more details and references. Given a set $\Th\subset\R^n$ with $\ox\in\Th$, the {\em tangent cone} to $\Th$ at $\ox$ is defined by
\begin{eqnarray}\label{2.5}
T_\Th(\ox):=\big\{w\in\R^n\;\big|\;\exists\,t_k{\downarrow}0,\;\;w^k\to w\;\;\mbox{ as }\;k\to\infty\;\;\mbox{with}\;\;\ox+t_kw^k\in \Th\big\}.
\end{eqnarray}
If $\Th$ is convex, the {\em normal cone} to $\Th$ at $\ox\in\Th$ in the sense of convex analysis   is
\begin{equation*}
N_\Th(\ox):=\big\{v\in\R^n\;\big|\;\la v,x-\ox\ra\le 0\;\;\mbox{for all}\;x\in \Th\big\}.
\end{equation*}
In the case where $\Th=\Q$, the second-order cone \eqref{icecream}, we get, respectively, the expressions
\begin{equation*}
T_\Q(y)=\begin{cases}
\mathbb{R}^{m+1}&\textrm{if }\;y\in\inter\Q,\\
\Q&\textrm{if }\;y=0,\\
\big\{y'\in\mathbb{R}^{m+1}\;\big|\;\langle\widetilde y,y'\rangle\le 0\big\}\quad&\textrm{if }\;y\in(\bd\Q)\setminus\{0\},
\end{cases}
\end{equation*}
\begin{equation}\label{ncone}
N_\Q(y)=\begin{cases}
\{0\}&\textrm{if }\;y\in\inter\Q,\\
\nQ&\textrm{if }\;y=0,\\
\mathbb{R}_+\widetilde y\quad&\textrm{if }\;y\in(\bd\Q)\setminus\{0\}.
\end{cases}
\end{equation}

Given further an extended-real-valued function $\ph\colon\R^n\to\oR:=(-\infty,\infty]$, its domain and epigraph are defined, respectively, by
$$
\dom\ph:=\big\{x\in\R^n\;\big|\;\ph(x)<\infty\big\}\quad\mbox{and}\quad\epi\ph:=\big\{(x,\mu)\in\R^{n+1}\;\big|\;\mu\ge\ph(x)\big\}.
$$
Considering next a set-valued mapping $F\colon\R^n\tto\R^m$ with its domain and graph given by
$$
\dom F:=\big\{x\in\R^n\;\big|\;F(x)\ne\emp\big\}\;\mbox{ and }\;\gph F:=\big\{(x,y)\in\R^n\times\R^m\;\big|\;y\in F(x)\big\},
$$
we define the {\em graphical derivative} of $F$ at $(\ox,\oy)\in\gph F$ via the tangent cone \eqref{2.5} by
\begin{equation}\label{gd}
DF(\ox,\oy)(u):=\big\{v\in\R^m\;\big|\;(u,v)\in T_{\scriptsize{\gph F}}(\ox,\oy)\big\},\quad u\in\R^n.
\end{equation}

A mapping $F\colon\R^n\tto\R^m$ is called {\em calm} at $(\ox,\oy)\in\gph F$ if there exist $\ell\ge 0$ and neighborhoods $U$ of $\ox$ and $V$ of $\oy$ for which
\begin{equation}\label{calm}
F(x)\cap V\subset F(\ox)+\ell\|x-\ox\|\B\;\mbox{ whenever }\;x\in U.
\end{equation}
It is said that $F$ has the {\em isolated calmness} property at $(\ox,\oy)\in\gph F$ if \eqref{calm} holds with the replacement of $F(\ox)$ by $\{\oy\}$ on the right-hand side therein. Furthermore, $F$ has the {\em robust isolated calmness} property at $(\ox,\oy)$ if
\begin{equation}\label{ricalm}
F(x)\cap V\subset\{\oy\}+\ell\|x-\ox\|\B\;\;\mbox{ with }\;F(x)\cap V\ne\emp\;\mbox{ for all }\;x\in U.
\end{equation}
Properties of this type go back to Robinson \cite{rob} who introduced the {\em upper Lipschitzian} version of calmness corresponding to \eqref{calm} with $V=\R^m$. Similarly to \eqref{ricalm}, we say that $F$ has the {\em robust isolated upper Lipschitzian} property if \eqref{ricalm} holds with $V=\R^m$. It is well known that \eqref{calm} is equivalent to the metric subregularity of the inverse mapping $F^{-1}$ at $(\oy,\ox)$. These ``one point" properties are more subtle and essentially less investigated than their robust ``two-points" counterparts (as metric regularity and Lipschitz-like/Aubin ones), while their importance for optimization theory, numerical algorithms, and applications has been broadly recognized in the literature; see, e.g., \cite{cst,dr,dl,g,ho,is14,m18,mms19b,yy} with the references and discussions therein.\vspace*{0.03in}

Turning now to the constructions of second-order variational analysis, for a function $\ph\colon\R^n\to\oR$, define the parametric
family of {\em second-order difference quotients} at $\ox$ for $\ov\in\R^n$   by
\begin{equation*}
\Delta_t^2\ph(\bar x,\ov)(w)=\dfrac{\ph(\ox+tw)-\ph(\ox)-t\langle\ov,w\rangle}{\frac{1}{2}t^2}\quad\mbox{with}\;\;w\in\R^{n},\;t>0.
\end{equation*}
If $\ph(\ox)$ is finite, the {\em second subderivative} of $\ph$ at $\ox$ for $\ov$ and $w$ is defined by
\begin{equation}\label{ssd}
\d^2\ph(\bar x,\ov)(w)=\liminf_{\substack{t\dn 0\\w'\to w}}\Delta_t^2\ph(\ox,\ov)(w').
\end{equation}

Following \cite[Definition~13.6]{RoW98}, a function $\ph\colon\R^n\to\oR$ is said to be {\em twice epi-differentiable} at $\bar x$ for $\ov$ if the sets $\epi\Delta_t^2\ph(\ox,\ov)$ converge to $\epi\d^2\ph(\bar x,\ov)$ as $t\downarrow 0$. If in addition the second subderivative is a proper function (i.e., does not take the value $-\infty$ and is finite at some point), then we say that $\ph$ is {\em properly} twice epi-differentiable at $\bar x$ for $\ov$. The twice epi-differentiability of $\ph$ at $\bar x$ for $\ov$  can be understood equivalently by \cite[Proposition~7.2]{RoW98} as that for every $w\in\R^n$ and every sequence $t_k\downarrow 0$ there exists a sequence $w^k\to w$ with
\begin{equation*}
\Delta_{t_k}^2\ph(\bar x,\ov)(w^k)\to\d^2\ph(\bar x,\ov)(w).
\end{equation*}
Twice epi-differentiability, together with a precise calculation of the second subderivative \eqref{ssd} of the augmented Lagrangian \eqref{aL} associated with \eqref{CP}, plays a major role in our developments. This property was introduced by Rockafellar in \cite{r88} who verified it for {\em fully amenable} compositions. Quite recently \cite{mms19a,MMS19,ms19}, the class of extended-real-valued functions satisfying this property has been dramatically enlarged by showing that twice epi-differentiability holds under {\em parabolic regularity}, which covers the SOCP setting; see more details in the cited papers.

When $\ph$ is the indicator function to the second-order cone $\Q$, it is shown in \cite[Theorem~3.1]{HMS18} that $\ph=\dd_\Q$ is properly twice epi-differentiable at any $\ox\in\Q$ for every $\ov\in N_\Q(\ox)$ and its second subderivative \eqref{ssd} is calculated by the precise formula
\begin{equation}\label{subq}
\d^2\delta_\Q(\bx,\ov)(w)=\begin{cases}
\delta_{K_\Q(\bx,\ov)}(w)&\textrm{if }\bx\in(\inter\Q)\cup\{0\},\\
\dfrac{\Vert\ov\Vert}{\Vert\bx\Vert}\Big(\Vert w_r\Vert^2-w_0^2\Big)+\delta_{K_\Q(\bx,\ov)}(w)&\textrm{if }\;\bx\in(\bd\Q)\setminus\{0\},
\end{cases}
\end{equation}
where $w=(w_0,w_r)\in\R\times\R^m$, and where $K_\Q(\bx,\ov):=T_\Q(\ox)\cap\{\ov\}^\bot$ stands for the {\em critical cone} to the set $\Q$ at $\ox\in\Q$ for any normal direction $\ov\in N_\Q(\ox)$.

The {\em Karush-Kuhn-Tucker} (KKT) optimality system associated with \eqref{CP} is given by
\begin{equation}\label{KKT}
\nabla_x L(x,\lambda)=\nabla f(x)+\nabla\Phi(x)^*\lambda=0,\quad\lambda\in N_\Q\big(\Phi(x)\big),
\end{equation}
where $L(x,\lambda):=f(x)+\langle\lambda,\Phi(x)\rangle$ is the (standard) {\em Lagrangian} of problem \eqref{CP} with $(x,\lm)\in\R^n\times\R^{m+1}$. For any $\ox\in\R^n$, define the set of {\em Lagrange multipliers} associated with $\ox$ by
\begin{equation}\label{ms}
\Lambda(\bx):=\left\{\lambda\in\R^{m+1}\;\big|\;\nabla_x L(\ox,\lambda)=0,\;\lambda\in N_\Q\big(\Phi(\bx)\big)\right\}.
\end{equation}

Our major attention to the second subderivative \eqref{ssd} in this paper is due to its ability to characterize the second-order growth condition in \eqref{CP} and thus to provide a second-order sufficient condition for strict local minimizers of this problem. To this end, we recall the corresponding result from \cite[Proposition~7.3]{MMS19} justifying such an application for SOCPs.\vspace*{-0.05in}

\begin{proposition}[\bf SOSC yields second-order growth]\label{prop2nd} Let $(\ox,\olm)\in\R^n\times\R^{m+1}$ be a solution to the
KKT system \eqref{KKT}, and let the second-order sufficient condition
\begin{equation}\label{SOSC}
\begin{cases}
\big\langle\nabla^2_{xx}L(\bx,\bl)w,w\big\rangle+\d^2\delta_\Q\big(\Phi(\bx),\olm\big)\big(\nabla\Phi(\ox)w\big)>0\\
\textrm{for all }\;w\in\R^n\setminus\{0\}\;\textrm{ with }\;\nabla\Phi(\bx)w\in K_\Q\big(\Phi(\bx),\bl\big)
\end{cases}
\end{equation}
hold. Then there exist positive numbers $\ell,\gamma$ such that the second-order growth condition
\begin{equation}\label{growth1}
f(x)\ge f(\bx)+\sm\ell\Vert x-\bx\Vert^2\quad\textrm{ for all }\;x\in\B_\gamma(\bx)\textrm{ with }\;\Phi(x)\in\Q
\end{equation}
is satisfied for the second-order cone program \eqref{CP}.
\end{proposition}\vspace*{-0.05in}

Observe that the presented SOSC \eqref{SOSC} is equivalent to the second-order conditions used for SOCPs in other publications \cite{br,HMS18,KaS17}. This indeed follows from the second subderivative formula \eqref{subq}. Note also  that SOSC \eqref{SOSC} is  stronger than the conventional second-order sufficient condition for \eqref{CP}, the latter requires   the supremum of the quadratic term in \eqref{SOSC} over all the Lagrange multipliers from \eqref{ms} be positive. This stronger condition is in fact equivalent to the second-order growth \eqref{growth1} under an appropriate constraint qualification; see \cite[Theorem~7.2]{MMS19}. Let us now provide an equivalent version of SOSC \eqref{SOSC} that is often used in what follows.\vspace*{-0.05in}

\begin{remark}[\bf equivalent version of SOSC]\label{rem1}{\rm It is not hard to check that the formulated SOSC \eqref{SOSC} amounts to saying that there exists a number $\bar\ell>0$ such that we have
\begin{equation}\label{SOSC2}
\begin{cases}
\big\langle\nabla^2_{xx}L(\bx,\bl)w,w\big\rangle+\d^2\delta_\Q\big(\Phi(\bx),\olm\big)\big(\nabla\Phi(\ox)w\big)\ge\bar\ell\,\|w\|^2\\
\textrm{for all }\;w\in\R^n\;\textrm{ with }\;\nabla\Phi(\bx)w\in K_\Q\big(\Phi(\bx),\bl\big).
\end{cases}
\end{equation}
Conversely, the fulfillment of \eqref{SOSC2} at $(\ox,\olm)$ ensures that for any $\ell\in(0,\bar\ell)$ there exists a positive number $\gamma$ such that the second-order growth condition \eqref{growth1} is satisfied at $\ox$.}
\end{remark}\vspace*{-0.05in}

We conclude this section by recalling some properties of the augmented Lagrangian \eqref{aL} that are used below; see, e.g., \cite[Exercise~11.56]{RoW98}.\vspace*{-0.05in}

\begin{proposition}[\bf properties of the augmented Lagrangian]\label{paug} For \eqref{aL} with $(x,\lm,\rho)\in\mathbb{R}^n\times\mathbb{R}^{m+1}\times(0,\infty)$ the following hold:\\[1ex]
{\bf(i)} The function $\rho\mapsto\L(x,\lambda,\rho)$ is nondecreasing.\\[1ex]
{\bf(ii)} The function $\lm\mapsto\L(x,\lambda,\rho)$ is concave.
\end{proposition}\vspace*{-0.05in}

It follows from the direct differentiation of \eqref{aL} that for any $\rho>0$ we have
\begin{eqnarray}\label{eq15}
\begin{array}{ll}
\nabla_x\L(x,\lambda,\rho)=\nabla f(x)+\nabla\Phi(x)^*\Pi_\nQ\big(\rho\Phi(x)+\lambda\big)),\\
\nabla_\lambda\L(x,\lambda,\rho)=\rho^{-1}\big[\Pi_\nQ\big(\rho\Phi(x)+\lambda\big)-\lambda\big],
\end{array}
\end{eqnarray}
which allows us to readily deduce that $(\ox,\olm)$ is a solution to the KKT system \eqref{KKT} if and only if for any $\rho>0$ this pair satisfies the equation
\begin{equation}\label{aLE}
\big(\nabla_x\L(x,\lambda,\rho),\nabla_\lm\L(x,\lambda,\rho)\big)=(0,0).
\end{equation}

Finally in this section, recall some properties of the projection mapping for the second-order cone $\Q$ that are extensively exploited throughout the paper:
\begin{itemize}[noitemsep,topsep=2pt]
\item[{\bf(P1)}]$p=\Pi_\Q(y)$ if and only if $p\in\Q$, $y-p\in\nQ$, and $\la y-p,p\ra=0$.
\item[{\bf(P2)}] For every $y\in\R^{m+1}$ we have $y=\Pi_\Q(y)+\Pi_\nQ(y)$.
\item[{\bf(P3)}] For every $y\in\R^{m+1}$ we have $\big\la\Pi_\Q(y),\Pi_\nQ(y)\big\ra=0$.
\item[{\bf(P4)}] $\lambda\in N_\Q(y)$ if and only if $\Pi_\Q(y+\lambda)=y$.
\end{itemize}\vspace*{-0.15in}

\section{Error Bounds for Perturbed KKT Systems of SOCPs}\label{sect03}\vspace*{-0.05in}

Here we derive an efficient error bound estimate for the KKT system of problem \eqref{CP} under the validity of SOSC \eqref{SOSC}. This is highly important for the subsequent results of the paper.

A crucial role of error bounds in convergence analysis of major numerical algorithms has been well understood in optimization theory; see, e.g., the books \cite{fp,is14}. To the best of our knowledge, the first error bound estimate for KKT systems of NLPs under the classical second-order sufficient condition alone was derived in Hager and Gowda \cite[Lemma~2]{hg} and then was improved by Izmailov \cite{iz05} who replaced the conventional SOSC with the weaker {\em  noncriticality} of Lagrange multipliers introduced therein. It has been recently observed by Mordukhovich and Sarabi \cite{MoS19} that similar results for nonpolyhedral conic programs require an additional assumption of the {\em calmness} of Lagrange multiplier mappings associated with canonically perturbed KKT systems. The latter assumption automatically holds for NLPs.

For any fixed $\ox\in\R^n$ the {\em multiplier mapping} $M_\bx\colon\mathbb{R}^n\times\mathbb{R}^{m+1}\rightrightarrows\mathbb{R}^{m+1}$, associated with the canonically perturbed KKT system \eqref{KKT} of \eqref{CP}, is defined by
\begin{equation}\label{mulmap}
M_\bx(v,w):=\big\{\lambda\in\R^{m+1}\;\big|\;\nabla_x L(\bx,\lambda)=v,\;\lambda\in N_\Q\big(\Phi(\bx)+w\big)\big\},\;\;(v,w)\in\R^n\times\R^{m+1}.
\end{equation}
It is easy to see that $M_\bx(0,0)$ reduces to the set of Lagrange multipliers $\Lambda(\bx)$ of the unperturbed system \eqref{ms}. Given a solution $(\ox,\olm)$ to the KKT system \eqref{KKT}, the calmness condition \eqref{calm} for $M_{\ox}$ at $((0,0),\bl)$ reads as the existence of positive constants $\tau$ and $\gamma$ such that
\begin{equation*}
M_\bx(v,w)\cap\B_\gamma(\bl)\subset\Lambda(\bx)+\tau\big(\Vert w\Vert+\Vert v\Vert\big)\B\;\;\mbox{whenever}\;\;(v,w)\in\B_\gamma(0,0).
\end{equation*}
This can be equivalently rewritten as the existence of $\tau,\gamma>0$ such that the estimate
\begin{equation}\label{calmm}
\dist\big(\lambda;\Lambda(\bx)\big)\le\tau\big(\Vert\nabla_x L(\bx,\lambda)\Vert+\dist\big(\Phi(\bx);N_\Q^{-1}(\lambda)\big)\big)
\end{equation}
holds for all $\lambda\in\B_\gamma(\bl)$. We can easily check that for (polyhedral) NLPs the calmness of the multiplier mapping follows automatically from the classical Hoffman lemma. Efficient conditions for the calmness of \eqref{mulmap} in the SOCP framework \eqref{CP} are presented at the end of this section.\vspace*{0.03in}

Now we are ready to derive the main result of this section ensuring the aforementioned error bound estimate. Define the {\em residual function}   $\sigma\colon\mathbb{R}^n\times\mathbb{R}^{m+1}\to\mathbb{R}$ of the KKT system \eqref{KKT} by
\begin{equation}\label{res}
\sigma(x,\lambda):=\Vert\nabla_x L(x,\lambda)\Vert+\Vert\Phi(x)-\Pi_\Q\big(\Phi(x)+\lambda\big)\Vert,\;\;(x,\lm)\in\R^n\times\R^{m+1}.
\end{equation}
It is easy to see that if $(\ox,\olm)$ is a solution to the KKT system \eqref{KKT}, then it follows from property (P4) of the projection mapping that $\sigma(\ox,\olm)=0$. Using this and the Lipschitz continuity of $\sigma$ with respect to both $x$ and $\lambda$ around $(\bx,\bl)$, we can find constants $\gamma_2>0$ and $\kappa_2\ge 0$ such that
\begin{equation}\label{eb1}
\sigma(x,\lambda)\le\kappa_2\big(\Vert x-\bx\Vert+\dist(\lambda;\Lambda(\bx))\big)\quad\textrm{ for all}\;\;(x,\lambda)\in\B_{\gamma_2}(\bx,\bl).
\end{equation}
Below we show that the opposite inequality in \eqref{eb1}, which is crucial for our subsequent developments of the ALM, can be achieved if in addition both SOSC \eqref{SOSC} and the calmness of the multiplier mapping are satisfied.
The provided proof, being strongly based on the geometry of the second-order cone \eqref{icecream}, is much simpler than the  one given recently in \cite[Theorem~5.9]{MoS19} for ${\cal C}^2$-cone reducible cone programs that is based on a highly involved reduction technique. \vspace*{-0.05in}

\begin{theorem}[\bf error bound for SOCPs under calmness and SOSC]\label{error} Let $(\ox,\olm)$ be a solution to the KKT system \eqref{KKT}, and let SOSC \eqref{SOSC} hold at $(\bx,\bl)$. If the multiplier mapping $M_\bx$ in \eqref{mulmap} is calm at $((0,0),\bl)$, then there exist constants $\gamma_1>0$ and $\kappa_1\ge 0$ such that
\begin{equation}\label{eb}
\Vert x-\bx\Vert+\dist\big(\lambda;\Lambda(\bx)\big)\le\kappa_1\,\sigma(x,\lambda)\quad\mbox{ for all }\;\;(x,\lambda)\in\B_{\gamma_1}(\bx,\bl),
\end{equation}
where the residual function $\sigma$ is taken from \eqref{res}.
\end{theorem}\vspace*{-0.15in}
\begin{proof}
Observe that if $x=\bx$ and $\lambda\in\Lambda(\bx)$, then \eqref{eb} holds since both sides are equal to $0$. Let us now verify \eqref{eb} while assuming that either $x\ne\bx$ or $\lambda\notin\Lambda(\bx)$. We first show that
\begin{equation}\label{ebb}
\Vert x-\bx\Vert=O\big(\sigma(x,\lambda)\big)\quad\textrm{as}\quad(x,\lambda)\to(\bx,\bl).
\end{equation}
Arguing by contradiction, suppose that there exists a sequence $(\xk,\lk)\to(\bx,\bl)$ with either $\xk\ne\bx$ or $\lk\notin\Lambda(\bx)$  satisfying the strict inequalities
\begin{equation*}
\Vert\xk-\bx\Vert>k\,\sigma(\xk,\lk)>0\quad\mbox{for all}\;\,k\in \N,
\end{equation*}
which imply that $\sigma(\xk,\lk)=o(\Vert\xk-\bx\Vert)$. By the definition of $\sigma$ the latter means that
\begin{equation}\label{eb3}
\nabla_x L(\xk,\lk)=o(\Vert\xk-\bx\Vert)\quad\textrm{and}\quad\alk:=\Phi\big(\xk)-\Pi_\Q(\Phi(\xk)+\lk\big)=o(\Vert\xk-\bx\Vert).
\end{equation}
Using the second equality in \eqref{eb3} combined with property (P1), we get the relationships
\begin{equation}\label{eb5}
\Phi(\xk)-\alk\in\Q,\quad\lk+\alk\in\nQ,\quad\textrm{and}\quad\big\langle\Phi(\xk)-\alk,\lk+\alk\big\rangle=0,
\end{equation}
which in turn bring us to the inclusion
\begin{equation}\label{eb2}
\lk+\alk\in N_\Q\big(\Phi(\xk)-\alk\big).
\end{equation}
It follows from the calmness estimate \eqref{calmm} that
\begin{equation*}\label{eb14}
\dist\big(\lk+\alk;\Lambda(\bx)\big)\le\tau\big(\Vert\nabla_x L(\bx,\lk+\alk)\Vert+\dist(\Phi(\bx);N^{-1}_\Q(\lk+\alk))\big)
\end{equation*}
for all $k\in\N$ sufficiently large. Since the gradient $\nabla f$ and Jacobian $\nabla\Phi$ mappings are Lipschitz continuous around $\ox$, we always have the estimate
\begin{eqnarray*}
\Vert\nabla_x L(\bx,\lk+\alk)\Vert&\le&\Vert\nabla f(\xk)-\nabla f(\bx)\Vert+\Vert\nabla_x L(\xk,\lk)\Vert+\Vert(\nabla\Phi(\xk)-\nabla\Phi(\bx))^*\lk\Vert\nonumber\\
&&+\Vert \nabla\Phi(\bx)^*\alk\Vert=O(\Vert\xk-\bx\Vert).
\end{eqnarray*}
On the other hand, it follows from \eqref{eb2} that $\Phi(\xk)-\alk\in N^{-1}_\Q(\lk+\alk)$, and hence
\begin{equation*}
\dist\big(\Phi(\bx);N^{-1}_\Q(\lk+\alk)\big)\le\Vert\Phi(\xk)-\alk-\Phi(\bx)\Vert=O(\Vert\xk-\bx\Vert),
\end{equation*}
where the last equality comes from the Lipschitz continuity of $\Phi$ around $\bx$ and the condition $\alk=o(\Vert\xk-\bx\Vert)$. This ensures in turn that $\lk-\Hat\lambda^k=O(\Vert\xk-\bx\Vert)$, where $\hlk:=\Pi_{\Lambda(\bx)}(\lk)$. Passing to subsequences if necessary gives us
\begin{equation}\label{eb4}
\dfrac{\xk-\bx}{\Vert\xk-\bx\Vert}\to\xi\ne 0\quad\textrm{and}\quad\dfrac{\lk-\hlk}{\Vert\xk-\bx\Vert}\to\eta\;\mbox{ as }\;k\to\infty.
\end{equation}
Appealing now to the first estimate in \eqref{eb3}, we arrive at the equalities
\begin{eqnarray*}
o(\Vert \xk-\bx\Vert )&=&\nabla_x L(\xk,\lk)=\nabla_x L(\xk,\bl)+\nabla\Phi(\xk)^*(\lk-\bl)\nonumber\\
&=&\nabla_x L(\bx, \bl)+\nabla^2_{xx}L(\bx,\bl)(\xk-\bx)+o(\Vert\xk-\bx\Vert)\nonumber\\
&&+ \nabla\Phi(\xk)^*(\lk-\hlk)+\big(\nabla\Phi(\xk)-\nabla\Phi(\bx)\big)^*(\hlk-\bl)\nonumber\\
&=&\nabla^2_{xx}L(\bx,\bl)(\xk-\bx)+\nabla\Phi(\xk)^*(\lk-\hlk)\nonumber\\
&&+\big(\nabla^2\Phi(\bx)(\xk-\bx)+o(\Vert\xk-\bx\Vert)\big)^*(\hlk-\bl)+o(\Vert\xk-\bx\Vert).
\end{eqnarray*}
Dividing both sides by $\Vert\xk-\bx\Vert$ and then passing to the limit as $k\to\infty$ show that
\begin{equation}\label{eb9}
0=\nabla^2_{xx}L(\bx,\bl)\xi+\nabla\Phi(\bx)^*\eta.
\end{equation}

Let us now verify the inclusion $\nabla\Phi(\bx)\xi\in K_\Q(\Phi(\bx),\bl)=T_\Q(\Phi(\bx))\cap\{\bl\}^\perp$. Indeed, using the first relation in \eqref{eb5} yields
\begin{equation*}\label{eb6}
\Q\ni\Phi(\xk)-\alk=\Phi(\ox)+\|\xk-\ox\|\Big[\nabla\Phi(\bx)\Big(\frac{\xk-\bx}{\|\xk-\bx\|}\Big)+\frac{o(\Vert \xk-\bx\Vert)}{{\|\xk-\bx\|}}\Big],
\end{equation*}
which tells us that $\nabla\Phi(\ox)\xi\in T_\Q(\Phi(\bx))$. Combining this with $\olm\in N_\Q(\Phi(\bx))$, we obtain $\la\olm,\nabla\Phi(\ox)\xi\ra\le 0$. To prove the equality therein, deduce from \eqref{eb2} that
$$
0\ge\big\la\lk+\alk,\Phi(\ox)-\Phi(\xk)+\alk\big\ra=-\Big\la\lk+\alk,\|\xk-\ox\|\Big[\nabla\Phi(\bx)\Big(\frac{\xk-\bx}{\|\xk-\bx\|}\Big)+\frac{o(\Vert
\xk-\bx\Vert)}{{\|\xk-\bx\|}}\Big]\Big\ra.
$$
Dividing both sides by $\Vert\xk-\bx\Vert$ and then passing to the limit as $k\to\infty$ verify that $\la\olm,\nabla\Phi(\ox)\xi\ra\ge 0$. Thus we get $\la\olm,\nabla\Phi(\ox)\xi\ra=0$ and hence arrive at $\nabla\Phi(\bx)\xi\in K_\Q(\Phi(\bx),\bl)$.

Our next step is to prove the following inequality involving the second subderivative \eqref{ssd}:
\begin{equation}\label{eb10}
\big\langle\nabla\Phi\big(\bx)\xi,\eta\big\rangle\ge\d^2\dd_\Q\big(\Phi(\bx),\olm\big)\big(\nabla\Phi(\bx)\xi\big).
\end{equation}
To proceed, remember that $\hlk\in N_\Q(\Phi(\bx))$. Using \eqref{eb2} and the monotonicity of the normal cone mapping to a convex set, we get
\begin{eqnarray*}
0&\le&\big\langle\Phi(\xk)-\Phi(\bx)-\alk,\lk-\hlk+\alk\big\rangle\\
&=&\big\langle\nabla\Phi(\bx)(\xk-\bx)+o(\Vert\xk-\bx\Vert),\lk-\hlk+o(\Vert\xk-\bx\Vert\big\rangle.
\end{eqnarray*}
Dividing both sides by $\Vert\xk-\bx\Vert^2$ and passing to the limit as $k\to\infty$ give us
\begin{equation*}
\big\langle\nabla\Phi(\bx)\xi,\eta\big\rangle\ge 0.
\end{equation*}
This   combined with \eqref{subq} verifies \eqref{eb10} if either $\Phi(\bx)=0$,  $\Phi(x)\in\inter\Q$, or $\bl=0$.

It remains to validate \eqref{eb10} in the case where $\Phi(\bx)\in(\bd\Q)\setminus\{0\}$ and $\bl\ne 0$. Then \eqref{eb2} and the normal cone representation \eqref{ncone} allow us to find $t_k\in\R_+$ and $\hat t_k\in\R_+$ such that $\lk+\alk=t_k(\widetilde\Phi(\xk)-\tilde\alpha^k)$ and $\hlk=\hat t_k\widetilde\Phi(\bx)$ for large $k\in\N$. We clearly have $\lim_{k\to\infty}t_k=\lim_{k\to\infty}\hat t_k={\Vert\bl\Vert}/{\Vert\Phi(\bx)\Vert}$. Passing to a subsequence if necessary, assume without loss of generality that either $t_k\ge\hat t_k$ or $t_k\le\hat t_k$ for all $k\in\N$. If the former holds, then
\begin{eqnarray*}
&&\big\langle\Phi(\xk)-\Phi(\bx)-\alk,\lk-\hlk+\alk\big\rangle=\big\langle\Phi(\xk)-\Phi(\bx)-\alk,t_k\widetilde\Phi(\xk)-\hat t_k\widetilde\Phi(\bx)-t_k\tilde\alpha^k\big\rangle\\
&=&\hat t_k\big\langle\Phi(\xk)-\Phi(\bx)-\alk,\widetilde\Phi(\xk)-\widetilde\Phi(\bx)-\tilde\alpha^k\big\rangle+(t_k-\hat t_k)\big\langle\Phi(\xk)-\Phi(\bx)-\alk,\widetilde\Phi(\xk)-\tilde\alpha^k\big\rangle\nonumber\\
&=&\hat t_k\big\langle\nabla\Phi(\bx)(\xk-\bx)+o(\Vert\xk-\bx\Vert),\nabla\widetilde\Phi(\bx)(\xk-\bx)+o(\Vert\xk-\bx\Vert)\big\rangle-(t_k-\hat t_k)\big\langle\Phi(\bx),\widetilde\Phi(\xk)-\tilde\alpha^k\big\rangle\nonumber\\
&\ge&\hat t_k\big\langle\nabla\Phi(\bx)(\xk-\bx)+o(\Vert\xk-\bx\Vert),\nabla\widetilde\Phi(\bx)(\xk-\bx)+o(\Vert\xk-\bx\Vert)\big\rangle,
\end{eqnarray*}
where the third equality comes from $\Phi(\xk)-\alk\in\bd\Q$ and the last inequality is due to $\Phi(\bx)\in\Q$ while $\widetilde\Phi(\xk)-\tilde\alpha^k\in\nQ$. If the latter holds,  a similar argument brings us to
\begin{eqnarray*}
&&\big\langle\Phi(\xk)-\Phi(\bx)-\alk,\lk-\hlk+\alk\big\rangle=\big\langle\Phi(\xk)-\Phi(\bx)-\alk,t_k\widetilde\Phi(\xk)-\hat t_k\widetilde\Phi(\bx)-t_k\tilde\alpha^k\big\rangle\\
&=&t_k\big\langle\Phi(\xk)-\Phi(\bx)-\alk,\widetilde\Phi(\xk)-\widetilde\Phi(\bx)-\tilde\alpha^k\big\rangle+(t_k-\hat t_k)\big\langle\Phi(\xk)-\Phi(\bx)-\alk,\widetilde\Phi(\bx)\big\rangle\nonumber\\
&=&t_k\big\langle\nabla\Phi(\bx)(\xk-\bx)+o(\Vert\xk-\bx\Vert),\nabla\widetilde\Phi(\bx)(\xk-\bx)+o(\Vert\xk-\bx\Vert)\big\rangle+(t_k-\hat t_k)\big\langle\Phi(\xk)-\alk,\widetilde\Phi(\bx)\big\rangle\nonumber\\
&\geq&t_k\big\langle\nabla\Phi(\bx)(\xk-\bx)+o(\Vert\xk-\bx\Vert),\nabla\widetilde\Phi(\bx)(\xk-\bx)+o(\Vert\xk-\bx\Vert)\big\rangle.
\end{eqnarray*}
Dividing these estimates by $\Vert\xk-\bx\Vert^2$ and passing to the limit as $k\to\infty$ result in
\begin{equation*}
\big\langle\nabla\Phi(\bx)\xi,\eta\big\rangle\ge\frac{\Vert\bl\Vert}{\Vert \Phi(\bx)\Vert}\big\langle \nabla\Phi(\bx)\xi,\nabla\widetilde\Phi(\bx)\xi\big\rangle=\d^2\dd_\Q(\Phi(\ox),\olm)(\nabla\Phi(\ox)\xi),
\end{equation*}
where the last equality is taken from \eqref{subq}. This fully justifies \eqref{eb10}.

Combining now \eqref{eb10} with \eqref{eb9} implies that
\begin{equation*}
\big\langle\xi,\nabla^2_{xx} L(\bx,\bl)\xi\big\rangle+\d^2\dd_\Q\big(\Phi(\bx),\olm\big)\big(\nabla\Phi(\bx)\xi\big)\le\big\langle\xi,\nabla^2_{xx}L(\bx, \bl)\xi\big\rangle+\big\langle\nabla\Phi(\bx)\xi,\eta\big\rangle=0,
\end{equation*}
which contradicts the second-order sufficient condition \eqref{SOSC} since $\Phi(\bx)\xi\in T_\Q(\Phi(\bx))$ and $\xi\ne 0$, and  thus verifies   estimate \eqref{ebb}.

To finish the proof of the claimed error bound \eqref{eb}, it remains to show that
\begin{equation}\label{dis}
\dist\big(\lambda;\Lambda(\bx)\big)=O\big(\sigma(x,\lambda)\big)\quad\textrm{as}\quad(x,\lambda)\to(\bx,\bl).
\end{equation}
To proceed, pick $(x,\lm)$ satisfying \eqref{ebb} and denote $y:=\Pi_\Q(\Phi(x)+\lambda)-\Phi(x)$. Thus we get $\lm-y\in N_\Q(\Phi(x)+y)$. Moreover, since  $(x,\lambda)\to(\bx,\bl)$, we get $y\to 0$. Combining the latter with \eqref{calmm} readily yields the relationships
\begin{eqnarray*}
\dist\big(\lambda-y;\Lambda(\bx)\big)&=&O\big(\Vert\nabla_x L(\bx,\lambda-y)\Vert+\dist(\Phi(\bx);N_\Q^{-1}(\lambda-y))\big)\\
&=&O\big(\Vert\nabla_x L(x,\lambda)\Vert+\|y\|+\|x-\ox\|\big)=O\big(\sigma(x,\lambda)\big),
\end{eqnarray*}
where the last equality comes from \eqref{ebb}. Since
$$
\dist\big(\lambda;\Lambda(\bx)\big)-\dist\big(\lambda-y;\Lambda(\bx)\big)=O(\|y\|)=O\big(\sigma(x,\lambda)\big),
$$
we arrive at \eqref{dis}. The error bound \eqref{eb} follows from the combination of \eqref{ebb} and \eqref{dis}, and hence  completes the proof of the theorem. \end{proof}\vspace*{-0.05in}

Next we present an example showing that the assumed calmness of the multiplier mapping in Theorem~\ref{error} is essential for the validity of the error bound \eqref{eb}. In fact, the following example demonstrates more: not only does the {\em primal-dual} error bound \eqref{eb} fail without the calmness assumption on \eqref{mulmap}, but even the {\em primal estimate} \eqref{ebb} is violated in the absence of calmness. This illustrates a striking difference between NLPs and nonpolyhedral SOCPs.\vspace*{-0.05in}

\begin{example}[\bf failure of error bound in the absence of calmness of multiplier mappings] Consider SOCP \eqref{CP} with the data
$f\colon\mathbb{R}^2\to\mathbb{R}$ and $\Phi\colon\mathbb{R}^2\to\mathbb{R}^3$ defined by
\begin{equation*}
f(x):=x_2^2\quad\textrm{and}\quad\Phi(x):=(-x_1^2+x_2,x_2,0)\quad\mbox{with}\;\;x=(x_1,x_2)\in\mathbb{R}^2.
\end{equation*}
Take $\bx:=(0,0)$ and observe that $\Phi(\bx)=0$ and that
\begin{equation*}
\nabla f(\bx)=\begin{pmatrix}
0\\
0\\
\end{pmatrix},\quad\nabla\Phi(\bx)^*=\begin{bmatrix}
0&0&0\\
1&1 &0\\
\end{bmatrix},\quad\Lambda(\bx)=\nQ\cap\left\{(1,1,0)\right\}^\perp=\mathbb{R}_+(-1,1,0).
\end{equation*}
Letting $\bl:=(-1,1,0)\in\Lambda(\bx)$, we conclude that the pair $(\ox,\olm)$ satisfies the KKT system \eqref{KKT}. It follows from the equality
\begin{equation*}
\nabla_{xx}^2 L(\bx,\bl)=\nabla^2 f(\bx)+\nabla^2\langle\bl,\Phi\rangle(\bx)=2I_2,
\end{equation*}
with $I_2$ standing for the $2\times 2$ identity matrix, that SOSC \eqref{SOSC} holds at $(\bx,\bl)$. To show now that the multiplier mapping $M_\bx$ from \eqref{mulmap} is not calm at $((0,0),\bl)$, select $\lk:=\big(-1,t_k,\sqrt{1-t_k^2}\big)$ with $t_k\uparrow 1$ as $k\to\infty$, which yields $\lk\to\bl$ as $k\to\infty$ and $\lk\in\nQ$ for all $k\in\N$. Direct calculations give us the expressions
\begin{equation*}
\dist^2\big(\lk;\Lambda(\bx)\big)=\Big\Vert\lk-\dfrac{\langle\lk,\bl\rangle}{\Vert\bl\Vert^2}\bl\Big\Vert^2=\dfrac{3-2t_k-t_k^2}{2}\;\mbox{ and}
\end{equation*}
\begin{equation*}
\Vert\nabla f(\bx)+\nabla\Phi(\bx)^*\lk\Vert^2=(t_k-1)^2,
\end{equation*}
which lead us to the limit calculations
\begin{equation*}
\lim_{k\to\infty}\dfrac{\dist^2\big(\lk;\Lambda(\bx)\big)}{\Vert\nabla f(\bx)+\nabla\Phi(\bx)^*\lk\Vert^2}=\dfrac{1}{2} \lim_{k\to\infty}\dfrac{3-2t_k-t_k^2}{(t_k-1)^2}=\infty.
\end{equation*}
This tells us that the multiplier mapping $M_\bx$ is not calm at $((0,0),\bl)$.

Next we check that the primal estimate \eqref{ebb} fails in this example. To proceed, take $\xk:=(0,\alpha_k)$ with $\alpha_k:=-(t_k-1)/{2}$ and observe that $(\xk,\lk)\to(\bx,\bl)$ as $k\to\infty$. This yields
\begin{equation}\label{eq55}
\nabla f(\xk)+\nabla\Phi(\xk)^*\lk=\begin{pmatrix}
0\\
2\alpha_k\\
\end{pmatrix}+\begin{pmatrix}
0\\
t_k-1\\
\end{pmatrix}=\begin{pmatrix}
0\\
0\\
\end{pmatrix}=o(\Vert\xk-\bx\Vert ).
\end{equation}
On the other hand, since $\Phi(\xk)$ is a nonzero point on the boundary of $\Q$ and $\lm^k$ is a nonzero point on the boundary of $\nQ$, it follows that
\begin{equation*}
\Phi(\xk)+\lk=\alpha_k(1,1,0)+\Big(-1,t_k,\sqrt{1-t_k^2}\Big)\notin\Q\cup\nQ.
\end{equation*}
Letting $y^k:=\Phi(\xk)+\lk$, we calculate that
\begin{equation*}
\Pi_\Q(y^k)=\dfrac{1}{2}\Big(y_0^k+\Vert y_r^k\Vert\Big)\Big(1,\dfrac{y_r^k}{\Vert y_r^k\Vert}\Big)
\end{equation*}
and then easily check as $k\to\infty$ that
\begin{equation*}
\bigg(1,\dfrac{y_r^k}{\Vert y_r^k\Vert}\bigg)\to(1,1,0)\quad\mbox{and}\quad
\lim_{k\to\infty}\dfrac{y_0^k+\Vert y_r^k\Vert}{\Vert\xk-\bx\Vert}=\lim_{k\to\infty}\dfrac{\alpha_k-1+\sqrt{\alpha_k^2+2\alpha_kt_k+1}}{\alpha_k}=2.
\end{equation*}
This allows us to compute the limits
\begin{eqnarray*}
\lim_{k\to\infty}\dfrac{\Vert\Phi(\xk)-\Pi_\Q(\Phi(\xk)+\lk)\Vert}{\Vert\xk-\bx\Vert}&=&
\lim_{k\to\infty}\Big\Vert\dfrac{\Phi(\xk)}{\alpha_k}-\dfrac{y_0^k+\Vert y_r^k\Vert}{2\Vert\xk-\bx\Vert}\Big(1,\dfrac{y_r^k}{\Vert y_r^k\Vert}\Big)\Big\Vert\\
&=&\Vert(1,1,0)-(1,1,0)\Vert=0.
\end{eqnarray*}
Combining the latter with \eqref{eq55} demonstrates that the primal estimate \eqref{ebb} and hence the error bound \eqref{eb} both fail in this simple example.
\end{example}\vspace*{-0.03in}

Let us now turn our attention to efficient conditions that ensure the fulfillment of the imposed calmness of the multiplier mapping \eqref{mulmap}. First we provide an improvement of a result established recently in \cite[Theorem~4.1]{MoS19}, which gives a complete characterization of the calmness property of \eqref{mulmap} together with the {\em uniqueness} of Lagrange multipliers in terms of the {\em dual qualification condition} that involves the graphical derivative \eqref{gd} of the normal cone mapping for \eqref{icecream}. To proceed, consider the {\em fully perturbed} set of Lagrange multipliers $M\colon\R^n\times\R^n\times\R^{m+1}\tto\R^{m+1}$, where---in contrast to $M_{\ox}(v,w)$ in \eqref{mulmap}---the decision variable $x$ is also included in the perturbation procedure. We define this mapping by
\begin{equation}\label{mulmapx}
M(x,v,w):=\big\{\lambda\in\R^{m+1}\;\big|\;\nabla_x L(x,\lambda)=v,\;\lambda\in N_\Q\big(\Phi(x)+w\big)\big\}
\end{equation}
for $(x,v,w)\in\R^n\times\R^n\times\R^{m+1}$ and observe that $M(\bx,0,0)=M_{\bx}(0,0)=\Lambda(\bx)$.
The next proposition provides a full characterization of  the {\em upper Lipschitzian} property of the fully perturbed multiplier mapping $M$ via the dual qualification condition, which plays a key role in the convergent analysis of the ALM for SOCP \eqref{CP}.\vspace*{-0.05in}

\begin{proposition}[\bf calmness and uniqueness of Lagrange multipliers]\label{dual} Let $(\ox,\olm)$ be a solution to the KKT system \eqref{KKT}. Then the following assertions are equivalent:\\[1ex]
{\bf(i)} The  multiplier mapping $M_{\ox}$ is calm at $((0,0),\olm)$, and $\Lm(\ox)=\{\olm\}$, i.e., the mapping $M_{\ox}$ has the isolated calmness property at $(\ox,\olm)$.\\[1ex]
{\bf(ii)} We have the dual qualification condition
\begin{equation}\label{duq}
DN_\Q\big(\Phi(\bx),\bl\big)(0)\cap\ker\nabla\Phi(\bx)^*=\{0\}.
\end{equation}
{\bf(iii)} There exist positive numbers $\gamma_3$ and $\kappa_3$ such that  the upper Lipschitzian estimate
\begin{equation}\label{pt1}
M(x,v,w)\subset\{\olm\}+\kappa_3(\Vert x-\bx\Vert+\Vert v\Vert+\Vert w\Vert)\B\;\mbox{ for all  }\;(x,v,w)\in\B_{\gamma_3}(\bx,0,0)
\end{equation}
holds for the fully perturbed multiplier mapping \eqref{mulmapx}.
\end{proposition}\vspace*{-0.15in}
\begin{proof} The equivalence between (i) and (ii) was established in \cite[Theorem~4.1]{MoS19}. Also it is not hard to see that (iii) implies (i) since $M(\bx,0,0)=\Lambda(\bx)$. Thus it remains to verify the last implication (ii)$\implies$(iii). Observe to this end due to \cite[Corollary~3.4]{HMS18} that
\begin{equation*}
DN_\Q\big(\Phi(\bx),\bl\big)(0)=N_{K_\Q (\Phi(\bx),\bl)}(0)=K_\Q\big(\Phi(\bx),\bl\big)^*=\big(T_\Q\big(\Phi(\bx)\big)\cap\{\bl\}^\perp\big)^*,
\end{equation*}
which in turn yields the inclusion
\begin{equation*}
N_\Q\big(\Phi(\bx)\big)=T_\Q\big(\Phi(\bx)\big)^*\subset DN_\Q\big(\Phi(\bx),\bl\big)(0).
\end{equation*}
Then the dual qualification \eqref{dual} ensures the fulfillment of the basic  constraint qualification
\begin{equation*}
N_\Q\big(\Phi(\bx)\big)\cap\ker\nabla\Phi(\bx)^*=\{0\},
\end{equation*}
which implies that the Lagrange multiplier sets $M(x,v,w)$ are uniformly bounded for all $(x,v,w)$ in some neighborhood $\mathcal{U}$ of the nominal triple $(\bx,0,0)$.

Having this in hand and arguing by contraposition, suppose on the contrary that the upper Lipschitzian property \eqref{pt1} fails. The equivalence between (i) and (ii) readily implies that  $M(\bx,0,0)=\Lambda(\bx)=\{\bl\}$. Thus it follows from the contraposition assumption that there exist sequences of $(\xk,v^k,w^k)\to(\bx,0,0)$ as $k\to\infty$ and of the corresponding multipliers $\lk\in M(\xk,v^k,w^k)$ satisfying the inequality
\begin{equation}\label{pt2}
\Vert\lk-\bl\Vert>k(\Vert\xk-\bx\Vert+\Vert v^k\Vert+\Vert w^k\Vert)\;\mbox{ whenever }\;k\in\N.
\end{equation}
Suppose without loss of generality that $(\xk,v^k,w^k)\in\mathcal{U}$ for all $k\in\mathbb{N}$. Hence the sequence $\{\lk\}$ is bounded, and so it has a limiting point $\hl$. Taking into account the robustness (closed graph property) of the normal cone mapping $N_\Q$ with respect to perturbations of the initial point, the continuity of the mappings $\Phi,\nabla f$, and $\nabla\Phi$ as well as the convergence $(\xk,v^k,w^k)\to(\bx,0,0)$, we arrive at $\hl\in\Lambda(\bx)=\{\olm\}$, which tells us that $\lk\to\bl$ as $k\to\infty$. Letting now $t_k:=\Vert \lk-\bl\Vert$ ensures that $t_k\downarrow 0$ and allows us to conclude by \eqref{pt2} that
\begin{equation}\label{pt3}
\xk-\bx=o(t_k),v^k=o(t_k),\;\mbox{ and }\;w^k=o(t_k)\;\mbox{ as }\;k\to\infty.
\end{equation}
Furthermore, the passage to a subsequence if necessary gives us a vector $\eta\in\mathbb{R}^{m+1}\setminus\{0\}$ such that $\dfrac{\lk-\bl}{t_k}\to\eta$. Recalling that $\lk\in M(\xk,v^k,w^k)$, we get
\begin{eqnarray*}
o(t_k)=\;v^k&=&\nabla f(\xk)+\nabla\Phi(\xk)^*\lk\\
&=&\nabla f(\xk)-\nabla f(\bx)+\nabla\Phi(\xk)^*\lk-\nabla\Phi(\bx)^*\bl\\
&=&\nabla f(\xk)-\nabla f(\bx)+\big(\nabla\Phi(x^k)-\nabla\Phi(\bx)\big)^*\lk+\nabla\Phi(\bx)^*(\lk-\bl)\\
&=&o(t_k)+\nabla\Phi(\bx)^*(\lk-\bl),
\end{eqnarray*}
where the verification of the last equality uses the Lipschitz continuity of $\nabla f$ and $\nabla\Phi$ around $\bx$, the boundedness of $\{\lk\}$, and the first estimate in \eqref{pt3}. Dividing both sides of the latter by $t_k$ and passing to the limit as $k\to\infty$ result in $\eta\in\ker\nabla\Phi(\bx)^*$. On the other hand, we have
\begin{equation*}
\left(\dfrac{\Phi(\xk)+w^k-\Phi(\bx)}{t_k},\dfrac{\lk-\bl}{t_k}\right )=\dfrac{\big(\Phi(\xk)+w^k,\lk\big)-\big(\Phi(\bx),\bl\big)}{t_k}\in \dfrac{\gph N_\Q-\big(\Phi(\bx),\bl\big)}{t_k},
\end{equation*}
which yields $(0,\eta)\in T_{\gph N_\Q}(\Phi(\bx),\bl)$ and hence verifies the condition
\begin{equation*}
\eta\in DN_\Q\big(\Phi(\bx,\bl)\big)(0)\cap\ker\nabla\Phi(\bx)^*.
\end{equation*}
Since $\eta\ne 0$, the latter contradicts \eqref{duq} and thus justifies the claimed estimate \eqref{pt1}.
\end{proof}

A different sufficient condition for the upper Lipschitzian property \eqref{pt1} was obtained in \cite[Proposition~4.47]{bs} by using  a condition called  the ``strict constraint qualification."
This condition is strictly more restrictive than the dual qualification \eqref{duq}, which---as shown in Proposition~\ref{dual}---is indeed equivalent to the upper Lipschitzian estimate in \eqref{pt1}.\vspace*{0.03in}

Our next goal is to provide a more detailed analysis of the calmness of the multiplier mapping for \eqref{CP} entirely via the given SOCP data at the fixed solution $(\ox,\olm)$ to the KKT system \eqref{KKT}. Consider all the possible cases. If $\Phi(\bx)\in\inter\Q$, then it follows from the normal cone representation \eqref{ncone} that $\Lambda(\bx)=\{0\}$ for the set of Lagrange multipliers in \eqref{ms}. Since $M_{\ox}(0,0)=\Lm(\ox)$ and since $M_{\ox}(u,v)=\{0\}$ whenever the pair $(u,v)$ is sufficiently close to $(0,0)$, we surely get the calmness of the multiplier mapping at $((0,0),\olm)$ with $\olm=0$ in this case. If further $\Phi(\bx)\in(\bd\Q)\setminus\{0\}$, then it follows from \eqref{ncone} that $\Lambda(\bx)$ is the intersection of two polyhedral convex sets. Employing the classical Hoffman lemma ensures that
\begin{equation*}
\dist\big(\lambda;\Lambda(\bx)\big)=O\big(\|\nabla_x L(\bx,\lambda)\|+\dist(\lambda;N_\Q(\Phi(\ox))\big)=O\big(\|\nabla_x L(\bx,\lambda)\|+ \dist(\Phi(\ox);N^{-1}_\Q(\lm)\big)
\end{equation*}
for all $\lm$ close enough to $\olm\in N_\Q(\Phi(\ox))$, where the last equality comes from the fact that the mapping $N_\Q$ is clearly calm at $(\Phi(\ox),\olm)$ in this case. This again verifies the calmness property of the multiplier mapping \eqref{mulmap} at $((0,0),\olm)$.

Considering further the remaining case where $\Phi(\bx)=0$, we deduce from \cite[Proposition~4.1]{HMS18} that the set of Lagrange multipliers $\Lambda(\bx)$ admits one of the following representations:\\[1ex]
{\bf(a)} The strict complementarity holds for $\Lambda(\bx)$, i.e., $\Lambda(\bx)$ contains an interior point of $\nQ$.\\[1ex]
{\bf(b)} $\Lambda(\bx)=\{0\}$.\\[1ex]
{\bf(c)} $\Lambda(\bx)=\{\olm\}$ and $\olm\in\bd(\nQ)\setminus\{0\}$.\\[1ex]
{\bf(d)} $\Lambda(\bx)=\mathbb{R}_+\olm$ and $\olm\in\bd(\nQ)\setminus\{0\}$.\vspace*{0.1in}

The next proposition describes the calmness of multipliers for \eqref{CP} when $\Phi(\ox)=0$.\vspace*{-0.05in}

\begin{proposition}[\bf calmness of SOCP multipliers at vertex]\label{soca} Let $(\ox,\olm)$ be a solution for the generalized KKT system \eqref{KKT}, and let $\Phi(\ox)=0$. The following hold:\\[1ex]
{\bf(i)} In cases {\rm (a)} and {\rm (b)} for $\Lambda(\ox)$ the multiplier mapping $M_{\ox}$ is calm at $((0,0),\olm)$.\\[1ex]
{\bf(ii)} In case {\rm(c)} for $\Lambda(\ox)$ the calmness of $M_{\ox}$ at $((0,0),\olm)$ is equivalent to the full rank of $\nabla\Phi(\ox)$.
\end{proposition}\vspace*{-0.15in}
\begin{proof}
In case (a) we get from \cite[Proposition~4.1]{HMS18} that estimate \eqref{calmm} is satisfied, which verifies the claimed calmness property of the multiplier mapping. In case (b) it follows from \eqref{ms} that $\nabla f(\bx)=0$, which yields the equalities
\begin{equation}\label{eq51}
\nQ\cap\ker\nabla\Phi(\bx)^*=N_\Q\big(\Phi(\bx)\big)\cap\ker\nabla\Phi(\bx)^*=\Lambda(\bx)=\{0\},
\end{equation}
and so $\bl=0$ and $K_\Q(\Phi(\bx),\bl)=T_\Q(0)=\Q$. By \cite[Corollary~3.4]{HMS18} we have
$$
DN_\Q\big(\Phi(\bx),\bl\big)(0)=N_{K_\Q(\Phi(\bx),\bl)}(0)=\nQ.
$$
This together with \eqref{eq51} tells us the dual qualification condition \eqref{duq} holds in this case. Employing Proposition~\ref{dual}  confirms the calmness of the multiplier mapping $M_\bx$ at $((0,0),\bl)$.

Finally, consider case (c). If $\nabla\Phi(\bx)$ has full rank, then the dual qualification condition \eqref{duq} is satisfied. Hence
Proposition~\ref{dual} ensures that the multiplier mapping $M_\bx$ is calm at $((0,0),\bl)$. Conversely, the validity of the calmness property for $M_{\ox}$ in the framework of (c) implies by Proposition~\ref{dual} that the dual qualification condition \eqref{duq} holds. Combining this with the fact that $\olm\in\bd(\nQ)\setminus\{0\}$ in (c) confirms that   the matrix $\nabla\Phi(\bx)$ has full rank; see \cite[Theorem~4.5]{HMS18} for the verification of this claim. This completes the proof of the proposition.
\end{proof}\vspace*{-0.05in}

The above discussions paint a clear picture for the calmness of the multiplier mapping in all the possible cases but (d). It has not been clarified at this stage how to provide verifiable conditions ensuring the calmness property of $M_{\ox}$ in case (d).\vspace*{-0.15in}

\section{Second-Order Variational Analysis of Augmented Lagrangians}\label{paL}\vspace*{-0.05in}

This section aims at providing characterizations of the second-order growth condition for the penalized problem \eqref{aup}. Our main device to obtain such characterizations is the {\em second subderivative}. As observed by Rockafellar \cite[Theorem~2.2]{r89}, the second-order growth condition for a proper extended-real-valued function can be characterized via its second subderivative. Using this rather simple albeit powerful result for the penalized problem \eqref{aup} requires the {\em calculation} of the second subderivative of the augmented  Lagrangian \eqref{aL}.

We begin with the following assertion that calculates the second subderivative of the Moreau envelope of a convex function. Given $\ph\colon\R^n\to\oR$ and $\rho>0$, recall that the {\em Moreau envelope} of $\ph$ relative to $\rho$ is defined by the infimal convolution
\begin{equation}\label{moreau}
(e_{1/\rho}\ph)(x):=\inf_w\Big\{\ph(w)+\sm\rho\|w-x\|^2\Big\},\quad x\in\R^n.
\end{equation}\vspace*{-0.25in}

\begin{proposition}[\bf second subderivatives of Moreau envelopes]\label{more} Let $\ph\colon\R^n\to\oR$ be a proper, lower semicontinuous {\rm(}l.s.c.{\rm)}, and convex function, and let $\ov\in\partial\ph(\ox)$. If $\ph$ is twice epi-differentiable at $\ox$ for $\ov$, then for any $\rho>0$ the Moreau envelope $e_{1/\rho}\ph$ is properly twice epi-differentiable at $\ox+\rho^{-1}\ov$ for $\ov$ and its second subderivative at this point is calculated by
\begin{equation}\label{mss}
\d^2 (e_{1/\rho}\ph)(\ox+\rho^{-1}\ov,\ov)(w)=e_{1/{2\rho}}\big(\d^2\ph(\ox,\ov)\big)(w)\;\;\mbox{for all }\;w\in\R^n.
\end{equation}
\end{proposition}\vspace*{-0.05in}
\begin{proof} Fix $\rho>0$. It follows from \cite[Theorem~11.23]{RoW98} that
\begin{equation}\label{eq31}
(e_{1/\rho}\ph)^*(z)=\ph^*(z)+\sm\rho^{-1}\Vert z\Vert^2\;\mbox{ for all }\;z\in\R^n,
\end{equation}
where `$*$' signifies the Fenchel conjugate in the sense of convex analysis. Because $\ph$ is proper, convex, and twice epi-differentiable at $\ox$ for $\ov$, we deduce from \cite[Proposition~13.20]{RoW98} that $\d^2\ph(\ox,\ov)$ is proper, l.s.c., and convex as well. Furthermore, it follows from \cite[Theorem~13.21]{RoW98} that the proper twice epi-differentiability of $\ph$ at $\ox$ for $\ov$ yields this property for the conjugate function $\ph^*$ at $\ov$ for $\ox$. Employing \cite[Proposition~12.19]{RoW98} tells us that the inclusion $\ov\in\partial\ph(\ox)$ ensures that $\nabla (e_{1/\rho}\ph)(\ox+\rho^{-1}\ov)=\ov$. Combining these facts with \eqref{eq31} and the sum rule for twice epi-differentiability from \cite[Exercise~13.18]{RoW98} implies that $(e_{1/\rho}\ph)^*$ is properly twice epi-differentiable at  $\ov$ for $\ox+\rho^{-1}\ov$ and that its second subderivative is given by
\begin{equation}\label{eq32}
\d^2 (e_{1/\rho}\ph)^*(\ov,\ox+\rho^{-1}\ov)(w)=\d^2\ph^*(\ov,\ox)(w)+\rho^{-1}\Vert w\Vert^2\;\;\mbox{for all}\;w\in\R^n.
\end{equation}
This together with \cite[Theorem~13.21]{RoW98} yields the proper twice epi-differentiability of $e_{1/\rho}\ph$ at $\ox+\rho^{-1}\ov$ for $\ov$. Thus the second subderivative of the latter function can be calculated by
\begin{eqnarray*}
\sm \d^2 (e_{1/\rho}\ph)(\ox+\rho^{-1}\ov,\ov)(w)& =&\Big(\sm\d^2(e_{1/\rho}\ph)^*(\ov,\ox+\rho^{-1}\ov)\Big)^*(w)\\
&=&\inf_{u\in\R^n}\Big\{\big(\sm\d^2\ph^*(\ov,\ox)\big)^*(u)+\sm \rho\Vert u-w\Vert^2\Big\}\\
&=&\inf_{u\in\R^n}\Big\{\sm\d^2\ph(\ox,\ov)(u)+\sm\rho\Vert u-w\Vert^2\Big\},
\end{eqnarray*}
where the first equality comes from  \cite[Theorem~13.21]{RoW98}, the second one is due to \eqref{eq32} and \cite[Proposition~14.1(i)]{BaC11}, and the last equality follows from \cite[Theorem~13.21]{RoW98}. This readily justifies the claimed formula for the second subderivative of $e_{1/\rho}\ph$ at $\ox+\rho^{-1}\ov$ for $\ov$.
\end{proof}\vspace*{-0.05in}

The second subderivative of the Moreau envelope for general prox-regular functions was established in \cite[Exercise~13.45]{RoW98}. However, there are several differences between the latter result and Proposition~\ref{more}. Firstly, the result of \cite{RoW98} was obtained for $\ov=0$ and $\rho>0$ sufficiently large. Our result does not demand neither of these requirements. Secondly, there is the coefficient $1/2$ in \cite[Exercise~13.45]{RoW98}, which does not appear in \eqref{mss}. The price for a nicer formula, however, is confining ourselves to the framework to convex functions.\vspace*{0.03in}

Proposition~\ref{more} allows us to obtain the required calculation of the second subderivative of the augmented Lagrangian \eqref{aL}.\vspace*{-0.05in}

\begin{theorem}[\bf second subderivatives of augmented Lagrangians]\label{ssc} Let $(\bx,\olm)$ be a solution to the KKT system \eqref{KKT}. Then for any $\rho>0$ the function $x\mapsto\L(x,\bl,\rho)$ defined via the augmented Lagrangian \eqref{aL} is twice epi-differentiable at $\bx$ for $0$ and its second subderivative is given by
\begin{equation}\label{eq40}
\d^2_x\L\big((\bx,\bl,\rho),0\big)(w)=\big\langle w,\nabla^2_{xx}L(\bx,\bl)w\big\rangle+Q_{\bx,\bl,\rho}(\nabla\Phi(\bx)w)+\rho\,\dist^2\big(\nabla\Phi(\bx)w;K_\Q(\Phi(\bx),\bl)\big),
\end{equation}
for  $w\in\R^n$, where the quadratic function $Q_{\bx,\bl,\rho}:\R^{m+1}\to\R$ is defined by
\begin{equation}\label{fq}
Q_{\bx,\bl,\rho}(v):=
\begin{cases}
0&\mbox{if}\;\Phi(\bx)\in(\inter\Q)\cup\{0\}\;\mbox{or}\;\bl=0,\\
\dfrac{\rho\Vert\bl\Vert}{\rho\Vert\Phi(\bx)\Vert+\Vert\bl\Vert}\Big(\Vert v_r\Vert^2-\dfrac{\langle\bl_r,v_r\rangle^2}{\Vert\bl_r\Vert^2}\Big)&\mbox{if}\;\Phi(\bx)\in(\bd\Q)\setminus\{0\}\;\mbox{and}\;\bl\ne 0.
\end{cases}
\end{equation}
\end{theorem}\vspace*{-0.15in}
\begin{proof} Since $(\bx,\olm)$ is a solution to the KKT system \eqref{KKT}, we have $\nabla_x\L(\ox,\bl,\rho)=0$, where $\nabla_x\L$ is calculated in \eqref{eq15}. The twice epi-differentiability of the function $x\mapsto\L(x,\bl,\rho)$ at $\ox$ for $\ov=0$ follows from \cite[Theorem~8.3(i)]{MMS19}. Let us proceed with the second subderivative calculation for the latter function. If either $\Phi(\bx)\in(\inter\Q)\cup\{0\}$ or $\bl=0$, then by \eqref{subq} we get
$$
\d^2\dd_\Q\big(\Phi(\ox),\olm\big)(w)=\dd_{K_\Q(\Phi(\ox),\olm)}(w)\;\mbox{ whenever }\;w\in\R^{m+1}.
$$
Employing again \cite[Theorem~8.3(i,iii)]{MMS19} and the second subderivative calculation \eqref{mss} from Proposition~\ref{more} for the Moreau envelope \eqref{moreau} of  $\ph=\dd_\Q$ tells us that
\begin{eqnarray*}
\d^2_x\L\big((\bx,\bl,\rho),0\big)(w)&=&\big\langle w,\nabla^2_{xx}L(\bx,\bl)w\big\rangle+e_{1/{2\rho}}\big(\d^2\dd_\Q(\Phi(\ox),\olm)\big)(w)\\
&=&\big\langle w,\nabla^2_{xx}L(\bx,\bl)w\big\rangle+\inf_{u\in\R^{m+1}}\big\{\dd_{K_\Q(\Phi(\ox),\olm)}(u)+\rho\|u-\nabla\Phi(\bx)w\|^2\big\}\\
&=&\big\langle w,\nabla^2_{xx}L(\bx,\bl)w\big\rangle+\rho\,\dist^2\big(\nabla\Phi(\bx)w;K_\Q(\Phi(\bx),\bl)\big),
\end{eqnarray*}
which verifies formula \eqref{eq40} with $Q_{\ox,\olm,\rho}(\nabla\Phi(\bx)w)=0$ from \eqref{fq} in this case. Assuming next that $\Phi(\bx)\in(\bd\Q)\setminus\{0\}$ and $\bl\ne 0$, define the function $\theta(y):=\sm\dist^2(y;\Q)$ for $y\in\R^{m+1}$. It is well known that $\theta$ is continuously differentiable on $\R^{m+1}$ and its gradient is given by
\begin{equation*}
\nabla\theta(y)=\Pi_\nQ(y)\;\;\mbox{whenever}\;\;y\in\R^{m+1}.
\end{equation*}
Since $\olm\in N_\Q(\Phi(\ox))$ and since $\Phi(\bx)\in(\bd\Q)\setminus\{0\}$ with $\bl\ne 0$, we get $\oy\notin\Q\cup(\nQ)$ with $\bar y=(\oy_0,\oy_r):=\Phi(\bx)+\rho^{-1}\bl$. This clearly yields $\|\oy_r\|>0$, and so we arrive at
\begin{equation*}
\nabla\theta(y)=\Pi_\nQ(y)=\frac{1}{2}\Big(1-\dfrac{y_0}{\Vert y_r\Vert}\Big)\big(-\Vert y_r\Vert,y_r\big)
=\frac{1}{2}\Big(y_0-\Vert y_r\Vert, y_r-y_0\dfrac{y_r}{\Vert y_r\Vert}\Big)
\end{equation*}
for all $y$ close to $\oy$. This confirms, in particular, that $\theta$ is $\mathcal{C}^2$-smooth around $\bar y$ with
\begin{equation}\label{eq36}
\nabla^2\theta(\oy)=\nabla\Pi_\nQ(\oy)=\dfrac{1}{2}\begin{pmatrix}
1&-\dfrac{\oy_r^*}{\Vert\oy_r\Vert}\\
-\dfrac{\oy_r}{\Vert\oy_r\Vert}&I_m-\dfrac{\oy_0}{\Vert\oy_r\Vert}I_m+\dfrac{\oy_0}{\Vert\oy_r\Vert}\dfrac{\oy_r\oy_r^*}{\Vert\oy_r\Vert^2}\\
\end{pmatrix},
\end{equation}
where $I_m$ the $m\times m$ identity matrix, and where $\oy_r^*$ stands for the corresponding vector row. Since $\Phi(\bx)\in(\bd Q)\setminus\{0\}$ and $\bl\in N_\Q(\Phi(\bx))\setminus\{0\}$, it follows that $\bl=t\widetilde\Phi(\bx)=t\left(-\Phi_0(\bx),\Phi_r(\bx)\right)$ for some $t>0$ and $\bl_0=-\Vert \bl_r\Vert $. Thus we have
\begin{equation*}
\bar y=\Phi(\bx)+\rho^{-1}\bl=\dfrac{1}{t}\big(-\bl_0,\bl_r\big)+\dfrac{1}{\rho}\big(\bl_0,\bl_r\big)=\bigg(\dfrac{t-\rho}{t\rho}\bl_0,
\dfrac{t+\rho}{t\rho}\bl_r\bigg)=\bigg(\dfrac{\rho-t}{t\rho}\Vert\bl_r\Vert,\dfrac{\rho+t}{t\rho}\bl_r\bigg).
\end{equation*}
Plugging the latter into \eqref{eq36} gives us the gradient formula
\begin{equation*}
\nabla\Pi_\nQ(\bar y)=\dfrac{1}{2}\begin{pmatrix}
1&-\dfrac{\bl_r^*}{\Vert\bl_r\Vert}\\
-\dfrac{\bl_r}{\Vert\bl_r\Vert}&\dfrac{2t}{\rho+t}I_m+\dfrac{\rho-t}{\rho+t}\dfrac{\bl_r\bl_r^*}{\Vert\bl_r\Vert^2}
\end{pmatrix},
\end{equation*}
which being combined with \eqref{eq36} and $\bl_0=-\Vert\bl_r\Vert$ results in
\begin{eqnarray}\label{eq39}
\big\la\nabla^2\theta(\oy)v,v\big\ra &=&
\dfrac{1}{2}\left(v_0^2-\dfrac{2v_0}{\Vert\bl_r\Vert}\langle\bl_r,v_r\rangle+\dfrac{2t}{\rho+t}\Vert v_r\Vert^2+\dfrac{\rho-t}{\rho+t}\dfrac{\langle\bl_r,v_r\rangle^2}{\Vert\bl_r\Vert^2}\right)\nonumber\\
&=&\dfrac{1}{2}\left[v_0^2-2v_0\dfrac{\langle\bl_r,v_r\rangle}{\Vert\bl_r\Vert}+\bigg(\dfrac{\langle\bl_r,v_r\rangle}{\Vert\bl_r\Vert}\bigg)^2 \right]+\dfrac{t}{\rho+t}\bigg(\Vert v_r\Vert^2-\dfrac{\langle\bl_r,v_r\rangle^2}{\Vert\bl_r\Vert^2}\bigg )\nonumber\\&=&\dfrac{(\bl_0 v_0)^2+2\hl_0v_0\langle\bl_r,v_r\rangle+\langle\bl_r,v_r\rangle^2}{2\Vert\bl_r\Vert^2}+\dfrac{t}{\rho+t}\bigg(\Vert v_r\Vert^2-\dfrac{\langle\bl_r,v_r\rangle^2}{\Vert\bl_r\Vert^2}\bigg)\nonumber\\
&=&\dfrac{\langle\bl,v\rangle^2}{\Vert\bl\Vert^2}+\dfrac{t}{\rho+t}\bigg(\Vert v_r\Vert^2-\dfrac{\langle\bl_r,v_r\rangle^2}{\Vert\bl_r\Vert^2}\bigg)\nonumber\\
&=&\dist^2\big(v;K_\Q(\Phi(\bx),\bl)\big)+\dfrac{\Vert\bl\Vert}{\rho\Vert\Phi(\bx)\Vert+\Vert\bl\Vert}\bigg(\Vert v_r\Vert^2-\dfrac{\langle\bl_r,v_r\rangle^2}{\Vert\bl_r\Vert^2}\bigg),
\end{eqnarray}
for all $v=(v_0, v_r)\in \R^{m+1}$. In the last equality we use the fact that $K_\Q(\Phi(\bx),\bl)=\{\bl\}^\perp$ and $\Vert\bl\Vert=t\Vert\Phi(\bx)\Vert$. It follows from the twice differentiability of $\theta$ at $\oy$ that the function $x\mapsto\L(x,\bl,\rho)$ is twice differentiable at $\ox$ with its second subderivative computed by
\begin{eqnarray*}
\d^2_x\L((\bx,\bl,\rho),0)(w)=\big\la\nabla^2_{xx}\L(\bx,\bl,\rho)w,w\big\ra=\big\langle w,\nabla^2_{xx}L(\bx,\bl)w\big\rangle+\rho\big\la\nabla^2\theta(\oy)v,v\big\ra
\end{eqnarray*}
with $v=\nabla\Phi(\bx)w$. Combining this and \eqref{eq39} gives us the claimed second subderivative formula in this case and thus finishes the proof of the theorem.
\end{proof}\vspace*{-0.05in}

Now we are ready to establish complete pointwise characterizations of the second-order growth condition for the penalized problem \eqref{aup} in terms of SOSC \eqref{SOSC} and the second subderivative of the augmented Lagrangian \eqref{aL}.\vspace*{-0.05in}

\begin{theorem}[\bf characterizations of second-order growth condition for augmented Lagrangians]\label{growth}
Let $(\ox,\olm)$ be a solution to the KKT system \eqref{KKT} for SOCP \eqref{CP}. Then the following assertions are equivalent:\\[1ex]
{\bf(i)} The second-order sufficient condition \eqref{SOSC} holds at $(\bx,\bl)$.\\[1ex]
{\bf(ii)} There exists a constant $\rhola>0$ such that for any $\rho\ge\rhola$ we have
\begin{equation}\label{eq26}
\d^2_x\L((\bx,\bl,\rho),0)(w)>0\;\mbox{ whenever }\;w\in\mathbb{R}^n\setminus\{0\}.
\end{equation}
{\bf(iii)} There exist positive constants $\rhola,\della$, and $\ellla$ such that for any $\rho\ge\rhola$ we have
\begin{equation}\label{eq16}
\L(x,\bl,\rho)\ge f(\bx)+\ellla\Vert x-\bar x\Vert^2\;\mbox{ for all }\;x\in\B_\della(\bx).
\end{equation}
\end{theorem}\vspace*{-0.05in}
\begin{proof} Since $(\ox,\olm)$ is a solution to the KKT system \eqref{KKT},  for all $\rho>0$ we have $\L(\bx,\bl,\rho)=f(\bx)$ and $\nabla_x\L(\bx,\bl,\rho)=0$. Assuming that (ii) holds, deduce from \cite[Theorem 13.24]{RoW98} that the second-order growth condition  \eqref{eq16} for $\rho=\rhola$ follows from \eqref{eq26} with the same constant $\rho$. Appealing now to Proposition~\ref{paug}(i) tells us that
$$
\L(x,\lambda,\rho)\ge\L(x,\lambda,\rhola)\;\mbox{ whenever }\;\rho\ge\rhola.
$$
This combined with \eqref{eq16} for $\rho=\rhola$ justifies the second-order growth condition for any $\rho\ge\rhola$ and thus verifies (iii). The opposite implication (iii)$\implies$(ii) follows directly from the definition of the second subderivative.\vspace*{0.05in}

Assume now that (ii) holds and let $\rho\ge\rhola$. To justify (i), pick $w\in\R^n\setminus\{0\}$ with $v:=\nabla\Phi(\ox)w\in K_\Q(\Phi(\ox),\olm)$. We next show that
\begin{equation}\label{mi12}
\d^2\delta_\Q\big(\Phi(\bx),\olm\big)\big(v\big) \geq Q_{\bx,\bl,\rho}(v)+\rho\,\dist^2\big(v;K_\Q(\Phi(\bx),\bl)\big)
\end{equation}
for all $\rho >0$. If either $\Phi(\bx)\in(\inter\Q)\cup\{0\}$ or $\bl=0$, we get from \eqref{subq} and \eqref{fq} that
\begin{equation} \label{mi14}
\d^2\delta_\Q\big(\Phi(\bx),\olm\big)\big(v\big) =Q_{\bx,\bl,\rho}(v)+\rho\,\dist^2\big(v;K_\Q(\Phi(\bx),\bl)\big) = 0,
\end{equation}
where the last equality stems from the fact that $\nabla\Phi(\ox)w\in K_\Q(\Phi(\ox),\olm)$. Otherwise, if $\Phi(\bx)\in(\bd\Q)\setminus\{0\}$ and $\bl\ne 0$, then we get that $K_\Q(\Phi(\bx),\bl)=\{\bl\}^\perp$. It follows from $v\in K_\Q(\Phi(\bx),\bl)$ and $\bl\in\bd(\nQ)\setminus\{0\}$ that
\begin{equation}\label{po1}
\langle v_{r}, \bl_{r}\rangle^2 = v_{0}^2\bl_{0}^2=v_{0}^2\|\bl_{r}\|^2\quad
\textrm{and}\quad \|v_{r}\|^2\geq v_{0}^2.
\end{equation}
We then deduce from \eqref{subq} and \eqref{fq} that
\begin{eqnarray}\label{mi15}
&&\d^2\delta_\Q\big(\Phi(\bx),\olm\big)\big(v\big) -Q_{\bx,\bl,\rho}(v)-\rho\,\dist^2\big(v;K_\Q(\Phi(\bx),\bl)\big) \nonumber\\
&=&\frac{\Vert \bl\Vert}{\Vert\Phi(\bx) \Vert}\big(\Vert v_{r}\Vert^2-v_{0}^2\big) - \frac{\rho\Vert\bl\Vert}{\rho\Vert\Phi(\bx)\Vert+\Vert\bl\Vert}\Big(\Vert v_r\Vert^2-\dfrac{\langle\bl_r,v_r\rangle^2}{\Vert\bl_r\Vert^2}\Big)\nonumber\\
&=&\dfrac{\Vert\bl\Vert^2}{\Vert\Phi(\bx)\Vert(\rho\Vert\Phi(\bx)\Vert+\Vert\bl\Vert)}\big(\Vert v_{r}\Vert^2-v_{0}^2\big)- \dfrac{\rho\Vert\bl\Vert}{\rho\Vert\Phi(\bx)\Vert+\Vert\bl\Vert}\Big(v_{0}^2-\dfrac{\langle\bl_r,v_{r}\rangle^2}{\Vert\bl_{r}\Vert^2}\Big)\nonumber\\
&\geq& 0,
\end{eqnarray}
where the last inequality is due to estimates in \eqref{po1}. Thus, we justify \eqref{mi12} for $v\in K_\Q(\Phi(\bx), \bl)$. Note that \eqref{mi12} is obvious if $v\notin K_\Q(\Phi(\bx), \bl)=\dom\d^2\delta_\Q\big(\Phi(\bx),\olm\big)$. Referring to \eqref{eq40}, we get SOSC \eqref{SOSC} from \eqref{eq26} and \eqref{mi12}. Thus we are done with (ii)$\implies$(i).\vspace*{0.05in}

To complete the proof of the theorem, it remains to verify implication (i)$\implies $(ii). Since the second subderivative is positive homogenous of degree 2, to prove \eqref{eq26} it is neccesary and sufficient to verify the condition: for all $\rho>0$ sufficiently large we get
\begin{equation}\label{eq43}
\d_x^2\L\big((\bx, \bl, \rho),0\big)(w)>0\quad\mbox{whenever}\;\;w\in\S.
\end{equation}
Assuming that (i) holds, we first justify the claim that \eqref{eq43} holds for all $w\in\S$ with $v:=\nabla\Phi(\bx)w\in K_\Q(\Phi(\bx),\bl)$. It is worth mentioning that the quadratic function (in $w$) on the left-hand side of SOSC \eqref{SOSC} must attain its minimum value on the compact set $\S$. Let $\ell_0$ denote such a value, then by \eqref{SOSC} we have $\ell_0>0$. We now show that
\begin{equation}\label{mi13}
Q_{\bx,\bl,\rho}(v)+\rho\,\dist^2\big(v;K_\Q(\Phi(\bx),\bl)\big) \geq \d^2\delta_\Q\big(\Phi(\bx),\olm\big)\big(v\big)-\frac{\ell_0}{2}
\end{equation}
for all $\rho>0$ sufficiently large. In the above proof of the implication (ii)$\implies$(i), it is proved that the latter holds for all $\rho>0$ whenever $\Phi(\bx)\in(\inter\Q)\cup\{0\}$ or $\bl=0$, see \eqref{mi14}. Turning now to the remaining case with $\Phi(\bx)\in(\bd\Q)\setminus\{0\}$ and $\bl\ne 0$. Recall from \eqref{po1} and \eqref{mi15} that
\begin{eqnarray*}
&&Q_{\bx,\bl,\rho}(v)+\rho\,\dist^2\big(v;K_\Q(\Phi(\bx),\bl)\big) - \d^2\delta_\Q\big(\Phi(\bx),\olm\big)\big(v\big)\nonumber\\
&=&-\dfrac{\Vert\bl\Vert^2}{\Vert\Phi(\bx)\Vert(\rho\Vert\Phi(\bx)\Vert+\Vert\bl\Vert)}\big(\Vert v_{r}\Vert^2-v_{0}^2\big)
\geq-\dfrac{\Vert\bl\Vert^2\|v\|^2}{\Vert\Phi(\bx)\Vert(\rho\Vert\Phi(\bx)\Vert+\Vert\bl\Vert)}\nonumber\\
&\geq&-\dfrac{\Vert\bl\Vert^2\|\nabla\Phi(\bx)\|^2}{\Vert\Phi(\bx)\Vert(\rho\Vert\Phi(\bx)\Vert+\Vert\bl\Vert)},
\end{eqnarray*}
where, in the last equality, we use the fact that $v=\nabla\Phi(\bx)w$ with $\|w\|=1$. Pick $\varrho_0>0$ such that the condition
\begin{equation}\label{po5}
\dfrac{\big(\Vert\nabla\Phi(\bx)\Vert\cdot\Vert\bl\Vert\big)^2}{\Vert\Phi(\bx)\Vert(\rho\Vert\Phi(\bx)\Vert+\Vert\bl\Vert)}\leq\frac{\ell_0}{2} \quad\textrm{for all }\rho\geq\varrho_0
\end{equation}
is fulfilled. Then \eqref{mi13} is satisfied for the case with $\Phi(\bx)\in(\bd\Q)\setminus\{0\}$ and $\bl\ne 0$, and therefore, for all possible position of $\Phi(\bx) \in \Q$ and $\bl \in N_\Q(\Phi(\bx))$ whenever $\rho\geq\varrho_0$. Referring to \eqref{eq40} and SOSC \eqref{SOSC}, we get by \eqref{mi13} that
\begin{equation}\label{mi16}
\d_x^2\L\big((\bx, \bl, \rho),0\big)(w) \geq \big\langle w, \nabla^2_{xx}L(\bx, \bl)w\big\rangle +\d^2\delta_\Q\big(\Phi(\bx),\olm\big)\big(v\big)-\frac{\ell_0}{2}\geq \frac{\ell_0}{2}
\end{equation}
for all $w\in \S$ with $v= \nabla\Phi(\bx)w \in K_\Q(\Phi(\bx),\bl)$ and for all $\rho\geq\varrho_0$, which just completes the verification of \eqref{eq43} for such $w$.\vspace*{0.05in}

Next we decompose the unit sphere into the two pieces:
\begin{equation*}
\S_+:=\big\{w\in\S\,\big|\;\big\langle w, \nabla^2_{xx}L(\bx, \bl)w\big\rangle +Q_{\bx,\bl,\varrho_0}(v)>0\big\}
\end{equation*}
and
\begin{equation*}
\S_-:=\big\{w\in\S\,\big|\;\big\langle w, \nabla^2_{xx}L(\bx, \bl)w\big\rangle +Q_{\bx,\bl,\varrho_0}(v)\le 0\big\},
\end{equation*}
where $\varrho_0$ is taken from \eqref{po5}. We see from \eqref{fq} that the function $\rho\mapsto Q_{\bx, \bl, \rho}(v)$ is nondecreasing on $\mathbb{R}_+$, then by \eqref{eq40} the estimate \eqref{eq43} is satisfied for any $w\in \S_+$  and any $\rho\geq\varrho_0$. Define the function $\vartheta\colon\S_-\to\R$ by
\begin{equation*}
\vartheta(w):=-\dfrac{\big\langle\nabla^2_{xx}L(\bx,\bl)w,w\big\rangle}{\dist^2\big(\nabla\Phi(\bx)w;K_\Q(\Phi(\bx),\bl)\big)},\;\;w\in\S_-.
\end{equation*}
Picking an arbitrary vector $w\in\S_-$, we conclude from the just proved claim that $\nabla\Phi(\bx)w\notin K_\Q(\Phi(\bx),\bl)$. This confirms that $\dist\big(\nabla\Phi(\bx)w;K_\Q(\Phi(\bx),\bl)\big)>0$. Also we get by \eqref{fq} that $Q_{\bx,\bl,\varrho_0}(v)$ is always positive, then $w\in \S_-$ implies that $\big\langle\nabla^2_{xx}L(\bx,\bl)w,w\big\rangle\le 0$. Thus  the function $\vartheta$ is continuous and nonnegative on the compact set $\S_-$, and hence its maximum value over this set, denoted by $\varrho_1$, is  finite and nonnegative. This demonstrates that for any $\rho>\varrho_1$ we have the estimate
\begin{equation*}
\big\langle\nabla^2_{xx}L(\bx,\bl)w,w\big\rangle+\rho\,\dist^2\big(\nabla\Phi(\bx)w;K_\Q(\Phi(\bx),\bl)\big)>0\;\mbox{ whenever }\;w\in\S_-.
\end{equation*}
This together  with the above estimate for the case of $w\in\S_+$ and $\rho>\varrho_0$ verifies \eqref{eq26} for all $w\in\R^n\setminus\{0\}$ and $\rho\geq \rhola>\max\{\varrho_0,\varrho_1\}$ and thus completes the proof of the theorem.
\end{proof}\vspace*{-0.05in}

Implication (i)$\implies$(iii) in Theorem~\ref{growth} was established by Rockafellar in \cite[Theorem~7.4]{r93} for nonlinear programming problems. His proof strongly exploits the geometry of NLPs and does not appeal to the second subderivative as in our proof. For the second-order cone programming problem \eqref{CP}, the aforementioned implication, not the established equivalencies in Theorem~\ref{growth}, was obtained in \cite[Proposition~10]{lz}, where in addition the strict complementarity and nondegeneracy conditions were imposed.\vspace*{0.05in}

To proceed further, observe that both constants $\ell_{\olm}$ and $\gamma_{\olm}$ in \eqref{eq16} depend on $\olm$. Now we are going to find additional assumptions that allow us to justify the second-order growth condition \eqref{eq16} for {\em all} $\lm\in\Lm(\ox)$ sufficiently close to $\olm$, where the aforementioned constants do not depend on $\lm$. This is crucial for the convergence analysis of the ALM in the case of {\em nonunique} Lagrange multipliers. The rest of this section is mainly focusing on achieving such a {\em uniform second-order growth condition} for the augmented Lagrangian \eqref{aL}.\vspace*{0.05in}

We begin with the following lemma, which provides a common constant $\ell_{\olm}$ that works for all $\lm$ sufficiently close to $\olm$. Then we derive a similar result for $\gamma_{\olm}$ in the proof of the next theorem.\vspace*{-0.05in}

\begin{lemma}[\bf uniform estimate for second subderivatives of augmented Lagrangians]\label{usog} Let $(\bx,\bl)$ be a solution to the KKT system \eqref{KKT}, and let SOSC \eqref{SOSC} hold at $(\ox,\olm)$. Then there exist positive constants $\rhola$, $\ell_1$, $\epsilon_0$ such that for all $\rho\ge\rhola$ and $\lm\in\Lambda(\bx)\cap\B_{\epsilon_0}(\bl)$ we have
\begin{equation}\label{eq44}
\d^2_x\L((\bx,\lm,\rho),0)(w)\ge\sm\ell_1\|w\|^2\;\mbox{ whenever }\;w\in\R^n.
\end{equation}
\end{lemma}\vspace*{-0.15in}
\begin{proof} Theorem~\ref{growth} gives us a constant $\rhola>0$ for which condition \eqref{eq26} holds when $\rho\ge\rhola$. Recall that the second subderivative is l.s.c.\ and positive homogenous of degree $2$. Owing to \eqref{eq40}, condition \eqref{eq26} amounts to the existence of a constant $\ell_1>0$ such that
\begin{equation}\label{eq41}
\d^2_x\L((\bx,\bl,\rhola),0)(w)=\big\langle w,\nabla^2_{xx}L(\bx,\bl)w\big\rangle+Q_{\bx,\bl,\rhola}(v)+\rhola \,\dist^2\big(v;K_\Q(\Phi(\bx),\bl)\big)\ge\ell_1
\end{equation}
for all $w$ from the unit sphere $\S\subset\R^n$ and $v=\nabla\Phi(\bx)w$, where the quadratic form $Q_{\bx,\bl,\rhola}(\cdot)$ is taken from \eqref{fq} with $\rho=\rhola$. Let us now verify the existence of $\epsilon_0>0$ so that for any $\lm\in\Lambda(\bx)\cap\B_{\epsilon_0}(\bl)$  we have
\begin{equation}\label{bu01}
\d^2_x\L((\bx,\lm,\rhola),0)(w)=\big\langle w,\nabla^2_{xx}L(\bx,\lambda)w\big\rangle+Q_{\bx,\lm,\rhola}(v)+\rhola\,\dist^2\big(v;K_\Q(\Phi(\bx),\bl)\big)\ge\frac{\ell_1}{2},\quad w\in\S,
\end{equation}
where $Q_{\bx,\lm,\rhola}(\cdot)$ is taken from \eqref{fq} with replacing $\bl$ by $\lm$ and $\rho$ by $\rhola$. We first observe that \linebreak $\big\langle w,\nabla^2_{xx}L(\bx,\lambda)w\big\rangle$ converges to $ \big\langle w,\nabla^2_{xx}L(\bx,\bl)w\big\rangle$ as $\lm\to\olm$ with $\lm\in \Lambda(\bx)$ uniformly for all $w\in\S$ due to the following estimate
\begin{equation*}
|\big\langle w,\nabla^2_{xx} L(\bx,\lm)w\big\rangle-\big\langle w,\nabla^2_{xx}L(\bx,\bl)w\big\rangle|\le\|\nabla^2\Phi(\ox)\|\cdot\|\lm-\olm\|.
\end{equation*}
We now prove the uniform convergence of $Q_{\bx,\lm,\rhola}(v)\to Q_{\bx,\bl,\rhola}(v)$ as $\lm{\to}\olm$ with $\lm\in \Lm(\ox)$ for all $v\in\nabla\Phi(\bx)(\S)$. It is obvious for the case with $\Phi(\ox)\in(\inter\Q)\cup\{0\}$, since quadratic forms reduce to $0$ by \eqref{fq}. Assume that $\Phi(\ox)\in(\bd\Q)\setminus\{0\}$. If $\olm=0$, $\lm\to\olm$ with $\lm\in \Lm(\ox)$, then it follows from \eqref{fq} that
\begin{eqnarray*}
|Q_{\bx,\lm,\rhola}(w)-Q_{\bx,\bl,\rhola}(w)|&=&\dfrac{\rhola\Vert\lm\Vert}{\rhola\Vert\Phi(\bx)\Vert+\Vert\lm\Vert}\Big(\Vert v_r\Vert^2-\dfrac{\langle\lm_r,v_r\rangle^2}{\Vert\lm_r\Vert^2}\Big)\\
&\le&\dfrac{\Vert\lm\Vert}{\Vert\Phi(\bx)\Vert}\Vert v\Vert^2
\le\dfrac{\Vert\nabla\Phi(\ox)\Vert^2}{\Vert\Phi(\bx)\Vert}\|\lm-\olm\|,
\end{eqnarray*}
which justifies the claimed uniform convergence in this case as well. Finally, assume that $\olm\ne 0$ and $\lm\to\olm$ with $\lm\in\Lm(\ox)$ and suppose without loss of generality that $\lm\ne 0$. Since $\lm, \olm\in\Lm(\ox)$ and $\Phi(\ox)\in(\bd \Q)\setminus\{0\}$, it follows from \eqref{ncone} that there exist positive constants $t$ and $\bar t$ such that $\lm=t\widetilde\Phi(\bx)$ and $\olm=\bar t\widetilde\Phi(\bx)$. These relationships result in the equality
\begin{equation*}\label{bu02}
\dfrac{\langle\lm_r,v_r\rangle^2}{\Vert\lm_r\Vert^2}=\dfrac{\langle\olm_r,v_r\rangle^2}{\Vert\olm_r\Vert^2}.
\end{equation*}
Using this together with \eqref{fq} brings us to the estimates
\begin{eqnarray*}
|Q_{\bx,\lm,\rhola}(v)-Q_{\bx,\bl,\rhola}(v)|&=&\Big|\dfrac{\rhola\Vert\lm\Vert}{\rhola\Vert\Phi(\bx)\Vert+\Vert\lm\Vert}-\dfrac{\rhola\Vert\olm\Vert}{\rhola\Vert\Phi(\bx)\Vert+\Vert\olm\Vert}\Big|\Big(\Vert v_r\Vert^2-\dfrac{\langle\olm_r,v_r\rangle^2}{\Vert\olm_r\Vert^2}\Big)\\
&\le&|t-\bar t|\|v\|^2=\frac{\|\lm -\bl\|}{\|\Phi(\bx)\|}\|v\|^2
\le\frac{\|\nabla\Phi(\bx)\|^2}{\Phi(\bx)}\|\lm-\olm\|,
\end{eqnarray*}
which again justify the claimed uniform convergence in this last case. Thus we find a number $\epsilon_1>0$ ensuring the uniform condition
\begin{equation}\label{eq45}
\big\langle w,\nabla^2_{xx} L(\bx,\lm)w\big\rangle+Q_{\bx,\lm,\rhola}(v)\geq\big\langle w,\nabla^2_{xx}L(\bx,\bl)w+Q_{\bx,\bl,\rhola}(v)-\frac{\ell_1}{4}
\end{equation}
whenever $w\in\S$ and $\lm\in\Lambda(\bx)\cap\B_{\epsilon_1}(\bl)$. Next we intend to verify the existence of $\epsilon_2>0$ such that
\begin{equation}\label{eq46}
\begin{cases}
\dist^2\big(\nabla\Phi(\bx)w;K_\Q(\Phi(\bx),\lm)\big)-\dist^2\big(\nabla\Phi(\bx)w;K_\Q(\Phi(\bx),\bl)\big)\ge-\dfrac{\ell_1}{4\rhola}\\
\textrm{for all }\;w\in\S\;\;\mbox{and all}\;\;\lm\in\Lambda(\bx)\cap\B_{\epsilon_2}(\bl).
\end{cases}
\end{equation}
To proceed, consider the following four possible locations of $\olm$ in $\nQ$:\\[1ex]
{\bf(a)} $\bl=0$. In this case we have
\begin{equation*}
K_\Q\big(\Phi(\bx),\bl\big)=T_\Q\big(\Phi(\bx)\big)\supset K_\Q\big(\Phi(\bx),\lm\big)
\end{equation*}
for all $\lm\in\Lm(\ox)$, which verifies the fulfillment of \eqref{eq46}.\\[1ex]
{\bf(b)} $\bl\in\inter(\nQ)$ with $\Phi(\bx)=0$. If $\lm$ is sufficiently close to $\bl$, then $\lm\in\inter(\nQ)$. This yields
\begin{equation*}
K_\Q\big(\Phi(\bx),\lm\big)=K_\Q\big(\Phi(\bx),\bl\big)=\{0\},
\end{equation*}
which immediately ensures that \eqref{eq46} holds.\\[1ex]
{\bf(c)} $\bl\in\bd(\nQ)\setminus\{0\}$ with $\Phi(\bx)\in\bd(\Q)\setminus\{0\}$. If $\lm\to\olm$ with $\lm\in\Lm(\ox)$, we get $\lm=t\bl$ for some $t>0$, which confirms  that
\begin{equation*}
K_\Q(\Phi(\bx),\lm)=K_\Q\big(\Phi(\bx),\bl\big).
\end{equation*}
This clearly justifies the claimed estimate \eqref{eq46}.\\[1ex]
{\bf(d)} $\bl\in\bd(\nQ)\setminus\{0\}$ with $\Phi(\bx)=0$. In this case, we have for all $\lm\in\Lambda(\bx)\setminus\{0\}$ that
\begin{equation*}
K_\Q\big(\Phi(\bx),\lm\big)=\begin{cases}
\mathbb{R}_+\widetilde\lm\quad&\textrm{if }\;\lm\in\bd(\nQ)\setminus\{0\},\\
\{0\}&\textrm{if }\;\lm\in\inter(\nQ),
\end{cases}
\end{equation*}
where the tilde-notation for the ice-cream cone is defined at the end of Section~1. Then \eqref{eq46} is obviously satisfied when $\lm\in\Lambda(\bx)\cap\inter(\nQ)$. Assume now that $\lm\in[\Lambda(\bx)\cap\bd(\nQ)]\setminus\{0\}$.  It is not hard to verify that for any such a $\lm$ we get
\begin{equation*}
\dist^2\big(v;K_\Q(\Phi(\bx),\lm)\big)=\Vert v\Vert^2-\dfrac{1}{\Vert\lm\Vert^2 }\big(\max\big\{0,\big\langle\widetilde\lm,v\big\rangle\big\}\big)^2.
\end{equation*}
It is worth mentioning that the function $\big(\max(0, t)\big)^2$ is $\mathcal{C}^1$ on the whole real line. It follows from the latter formula that $\lambda \in [\Lambda(\bx)\cap\bd(\nQ)]\setminus\{0\} \mapsto \dist^2\big(v;K_\Q(\Phi(\bx),\lm)\big)$ is a $\mathcal{C}^1$ function relative to the set $[\Lambda(\bx)\cap\bd(\nQ)]\setminus\{0\}$. Taking this into account and choosing $\lm$ to be sufficiently close to $\olm$ ensure the existence of $\epsilon_2>0$ for which the uniform estimate \eqref{eq46} is guaranteed. This completes the justification of \eqref{eq46} for all the possible cases.\vspace*{0.05in}

Finally, denote $\epsilon_0:=\min\{\epsilon_1,\epsilon_2\}$ with $\epsilon_1$ and $\epsilon_2 $ taken from \eqref{eq45} and \eqref{eq46}, respectively. Combining  \eqref{eq41}, \eqref{eq45}, and \eqref{eq46} tells us that estimate \eqref{bu01} is satisfied for any $\lm\in\Lambda(\bx)\cap\B_{\epsilon_0}(\bl)$. Thus for any such a multiplier $\lm$ we have
$$
\d^2_x\L\big((\bx,\lm,\rhola),0\big)(w)\ge\sm\ell_1\|w\|^2\;\textrm{ whenever }\;w\in\R^n.
$$
This together with \eqref{eq40} and the fact that $\rho\mapsto Q_{\bx, \bl, \rho}(v)$ is nondecreasing on $\R_+$ implies for any $\lm\in\Lambda(\bx)\cap\B_{\epsilon_0}(\bl)$ that
$$
\d^2_x\L\big((\bx,\lm,\rho),0\big)(w)\ge\sm\ell_1\|w\|^2\;\textrm{ for all }\;w\in\R^n\;\;\mbox{and all}\;\;\rho\ge\rhola,
$$
which therefore completes the proof of the lemma.
\end{proof}\vspace*{-0.05in}

Now we are ready to derive a uniform version of the second-order growth condition for \eqref{aL}.\vspace*{-0.05in}

\begin{theorem}[\bf uniform second-order growth condition for augmented Lagrangians]\label{ugrowth}
Let $(\bx,\bl)$ be a solution to the KKT system \eqref{KKT}, and let SOSC \eqref{SOSC} hold at $(\ox,\olm)$. Assume in addition that the Lagrange multiplier set $\Lambda(\bx)$ in \eqref{ms} is either a polyhedron, or that the multiplier $\bl$ belongs to the interior of $\nQ$. Then there are positive constants $\rhola,\della,\ela,\ellla$ such that for all $\lm\in\Lambda(\bx)\cap\B_\ela(\bl)$ and $\rho\ge\rhola$ we have the uniform second-order growth condition
\begin{equation}\label{eq14}
\L(x,\lm,\rho)\ge f(\bx)+\ellla\Vert x-\bar x\Vert^2\;\mbox{ whenever }\;x\in\B_\della(\bx).
\end{equation}
\end{theorem}\vspace*{-0.15in}
\begin{proof} Take the positive constants $\ell_1$, $\epsilon_0$, and $\rhola$ from Lemma~\ref{usog} for which \eqref{eq44} holds whenever $\lm\in\Lambda(\bx)\cap\B_{\epsilon_0}(\bl)$ and $\rho\ge\rhola$. Using \cite[Theorem~13.24]{RoW98} and remembering that $\L(\ox,\lm,\rhola)=f(\ox)$ for all $\lm\in\Lm(\ox)$, we deduce from \eqref{eq44} that for any $\lm\in\Lambda(\bx)\cap\B_{\epsilon_0}(\bl)$ there exists $\gamma_{\lm}>0$ ensuring the estimate
\begin{equation}\label{eq22}
\L(x,\lm,\rhola)\ge f(\bx)+\dfrac{\ell_1}{4}\Vert x-\ox\Vert^2\;\textrm{ whenever }\;x\in\B_{\gamma_\lm}(\ox),
\end{equation}
where the constant $\ell_1/4$ can be chosen the same for all the multipliers $\lm\in\Lambda(\bx)\cap\B_{\epsilon_0}(\bl)$. This comes from \eqref{eq44} and the proof of \cite[Theorem~13.24]{RoW98}; see also Remark~\ref{rem1} for a similar discussion. However, the radii of the balls centered at $\ox$ in \eqref{eq22} depend on $\lm$. It is shown below that we can find a common radius for all the multipliers $\lm\in\Lambda(\bx)$ that are sufficiently close to $\bl$. To proceed, define the function $\ph\colon\R^{m+1}\to\oR$ by
\begin{equation}\label{eq23}
\varphi(\lambda):=\sup_{x\in\B_\della(\ox)}\dfrac{f(\bx)-\L(x,\lambda,\rhola)}{\Vert x-\ox\Vert^2}+\delta_{\Lambda(\bx)\cap\B_{\epsilon_0}(\bl)}(\lambda),\;\;\lm\in\R^{m+1}.
\end{equation}
Proposition~\ref{paug}(ii) tells us that the function $\lambda\mapsto\L(x,\lm,\rhola)$ is concave. This together with the convexity of the set
$\Lambda(\bx)\cap\B_{\epsilon_0}(\bl)$ ensures that $\varphi$ in \eqref{eq23} is a convex function. Let us now verify that for any $\lm\in\Lambda(\bx)\cap\B_{\epsilon_0}(\bl)$ the value $\ph(\lm)$ is finite. To this end, pick such a multiplier $\lm$ and observe that if ${\gamma_\lm}\ge\della$ we get by \eqref{eq22} the estimates
\begin{equation*}
\varphi(\lm)\le\sup_{x\in\B_{{\gamma_\lm}}(\ox)}\dfrac{f(\bx)-\L(x,\lm,\rhola)}{\Vert x-\ox\Vert^2}\le-\dfrac{\ell_1}{4}.
\end{equation*}
In particular, this implies that $\varphi(\bl)\le-\ell_1/4$. If ${\gamma_\lm}<\della$, then
\begin{equation*}
\varphi(\lm)\le\max\left\{\sup_{x\in\B_{\gamma_\lm}(\ox)}\dfrac{f(\bx)-\L(x,\lm,\rhola)}{\Vert x-\ox\Vert^2},\max_{{\gamma_\lm}\le\Vert x-\ox\Vert\le\della}\dfrac{f(\bx)-\L(x,\lm,\rhola)}{\Vert x-\ox\Vert^2}\right\}<\infty,
\end{equation*}
where the first term inside the maximum does not exceed $-{\ell_1}/{4}$ because of \eqref{eq22}, and where the second term is finite since it is the maximum of a continuous function over a compact set. This implies that $\ph(\lm)$ is finite for all $\lm\in\Lambda(\bx)\cap\B_{\epsilon_0}(\bl)$, which ensures   that
$$
\dom\varphi=\Lambda(\bx)\cap\B_{\epsilon_0}(\bl).
$$
If $\bl\in\inter(-\Q)$, we get $\olm\in\ri\Lm(\ox)$, which clearly implies that $\olm\in\ri(\dom\ph)$. Since $\ph$ is convex, it is continuous at $\olm$ relative to its domain. Hence we find $\ela\in(0,\epsilon_0]$ such that
\begin{equation}\label{eq24}
\varphi(\lm)\le\varphi(\bl)+\dfrac{\ell_1}{8}\le-\dfrac{\ell_1}{8}\;\textrm{ for all }\;\lm\in\dom\ph\cap\B_{\ela}(\bl)=\Lm(\ox)\cap\B_{\ela}(\bl).
\end{equation}

Next we proceed to achieve a similar result when $\Lambda(\bx)$ is a polyhedral convex set. In this case the collection of Lagrange multipliers $\Lm(\ox)$ is either a ray on the boundary of $\nQ$, or a singleton. If the latter holds, we obtain $\Lm(\ox)=\{\olm\}$, and hence the uniform growth condition \eqref{eq14} follows directly from \eqref{eq16}. If $\Lm(\ox)$ is a ray on the boundary of $\nQ$, then  $\Lambda(\bx)\cap\B_{\epsilon_0}(\bl)$ is a segment. If now $\olm\ne 0$, then we get $\olm\in\ri[\Lambda(\bx)\cap\B_{\epsilon_0}(\bl)]=
\ri(\dom\ph)$. Arguing as above leads us to \eqref{eq24} in this case. Otherwise, $\olm$ is an endpoint of the aforementioned  segment,
and thus $\olm=0$. Let $\lm_e$ be the other endpoint. If $\varphi(\lm_e)\le\varphi(\olm)+{\ell_1}/{8}$, then \eqref{eq24} holds for $\ela:=\epsilon_0$, which follows from the convexity of $\varphi$. Otherwise, we have that $\varphi(\lm_e)>\varphi(\olm)+{\ell_1}/{8}$.  Denote
\begin{equation*}
\bar t:=\dfrac{\ell_1}{8\big(\varphi(\lm_e)-\varphi(\bl)\big)}\in(0,1)\;\textrm{ and }\;\lm_{\bar t}:=(1-\bar t)\bl+\bar t\lm_e.
\end{equation*}
Then using the convexity of $\ph$ tells us that
\begin{equation*}
\varphi(\lm_{\bar t})\le(1-{\bar t})\varphi(\bl)+\bar t\varphi(\lm_e)=\varphi(\bl)+\bar t\big(\varphi(\lm_e)-\varphi(\bl)\big)=\varphi(\bl)+\dfrac{\ell_1}{8}\le-\dfrac{\ell_1}{8},
\end{equation*}
which readily yields \eqref{eq24} with $\ela:=\Vert\lm_{\bar t}-\bl\Vert\in(0,\epsilon_0]$.
This  completes the verification of \eqref{eq24} with some constant $\ela\in(0,\epsilon_0]$ if either $\Lambda(\bx)$ is a polyhedral convex set, or $\bl\in \inter(\nQ)$. Consequently, it follows from \eqref{eq23} and  $\eqref{eq24}$ that
\begin{equation}\label{eq25}
\L(x,\lm,\rhola)\ge f(\bx)+\dfrac{\ell_1}{8}\Vert x-\bar x\Vert^2\;\textrm{ for all }\;x\in\B_\della(\bx)\;\textrm{ and }\;\lm\in\Lambda(\bx)\cap\B_{\ela}(\bl).
\end{equation}
Employing now Proposition~\ref{paug}(i) gives us the inequality
\begin{equation*}
\L(x,\lm,\rho)\ge\L(x,\lm,\rhola)\;\mbox{ for all}\;\;\rho\ge\rhola.
\end{equation*}
Combining this with \eqref{eq25} and setting $\ellla:={\ell_1}/{8}$ verify the uniform growth condition \eqref{eq14}.
\end{proof}\vspace*{-0.05in}

A similar result to Theorem~\ref{ugrowth} was derived in \cite[Proposition~3.1]{FeS12} for NLPs.   We are not familiar with any previous results on the uniform second-order growth condition \eqref{eq14} for SOCPs. As shown in the next section, the second-order growth conditions obtained above are crucial for developing the augmented Lagrangian method for this class of optimization problems.\vspace*{-0.15in}

\section{Well-Posedness and Convergence Analysis of ALM for SOCPs}\label{sect05}\vspace*{-0.05in}

In this concluding section of the paper we apply the suggested approach and results of second-order variational analysis (which are undoubtedly of their independent interest) to the convergence analysis of the {\em augmented Lagrangian method} for solving SOCPs \eqref{CP}.

The principal idea of the ALM for \eqref{CP} is to solve a sequence of {\em unconstrained minimization} problems for which the objective functions, at each iteration, are approximations of the augmented Lagrangian \eqref{aL}. Namely, given the current iteration $(x^k,\lm^k,\rho_k)$, the ALM solves the following unconstrained problem (called a {\em subproblem}):
\begin{equation}\label{ag1}
\mbox{minimize }\;\L(x,\lk,\rho_k)\;\mbox{ for }\;x\in\R^n
\end{equation}
for next primal iterate $\xkk$ and then use it to construct the next dual iterate $\lkk$. More specifically, we aim at solving the {\em stationary equation}
\begin{equation}\label{aLEs}
\nabla_x\L(x,\lk,\rho_k)=0
\end{equation}
for $\xkk$ and then to update the corresponding multiplier by $\lkk:=\Pi_\nQ(\rho_k\Phi(\xkk)+\lk)$.

Since solving \eqref{aLEs} is not easy in practice, it is more convenient to choose an approximate solution $\xkk$ satisfying the {\em approximate stationary condition}
\begin{equation}\label{xkk}
\Vert\nabla_x\L(\xkk,\lk,\rok)\Vert\le\ek
\end{equation}
with a given accuracy/tolerance $\ek\ge 0$. Following the conventional terminology of nonlinear programming, we say that the ALM is {\em exact} of $\ek=0$, i.e., the exact stationary equation \eqref{aLEs} is used, and {\em inexact} if \eqref{xkk} with $\ek>0$ is under consideration. In this paper we deal with both exact and inexact versions of the ALM by choosing an arbitrary accuracy $\ek\ge 0$ sufficiently small. The ALM   for \eqref{CP} is described as follows.\vspace*{-0.05in}

\begin{algorithm}[\bf augmented Lagrangian method for SOCPs]\label{Alg1}
Choose $(x^0,\lambda^0)\in\mathbb{R}^n\times \R^{m+1}$ and $\bar\rho>0$. Pick $\ek\to 0$ as $k\to\infty$ and $\rok$ with $\rok\ge \bar\rho$ for all $k$ and set $k:=0$. Then:\\[1ex]
{\bf(1)} If $(\xk,\lk)$ satisfies a suitable termination criterion, stop.\\[1ex]
{\bf(2)} Otherwise, find $\xkk$ satisfying \eqref{xkk} and update the Lagrange multiplier by
\begin{equation}\label{lkk}
\lkk:=\Pi_\nQ\big(\rok\Phi(\xkk)+\lk\big).
\end{equation}
{\bf(3)} Set $k\leftarrow k+1$ and go to Step~1.
\end{algorithm}\vspace*{-0.03in}

To perform the well-posedness and convergence analysis of Algorithm~\ref{Alg1}, we need to make sure first of all that the ALM is {\em well-defined}, i.e., its subproblems constructed in \eqref{ag1} are {\em solvable}. The following theorem reveals that the optimal solution mappings to subproblems \eqref{ag1} enjoy the robust isolated upper Lipschitzian property uniformly in $\rho$. This confirms, in particular, that subproblems \eqref{ag1} always admit a local optimal solution. Note that the developed proof of the theorem requires only the second-order growth condition \eqref{eq16}, which is based on SOSC \eqref{SOSC}, without any additional assumptions.\vspace*{-0.05in}

\begin{theorem}[\bf solvability and robust stability of subproblems in ALM]\label{solv} Let $\rhola$, $\della$, and $\ellla$ be positive constants for which the second-order growth condition \eqref{eq16} holds whenever $\rho\ge\rhola$. Then there exist constants $\ell >0$, $\hat\gamma\in(0,\gamma_{\olm}]$, and $\epsilon>0$ such that the local optimal solution mapping $S_\rho\colon\R^{m+1}\to\R^n$ defined by
\begin{equation}\label{eq177}
S_\rho(\lambda):=\argmin\big\{\L(x,\lambda,\rho)\;\big|\;x\in\B_{\hat\gamma}(\bx)\big\},\quad\lm\in\R^{m+1},
\end{equation}
satisfies, for all $\lm\in\B_\epsilon(\olm)$ and all $\rho\in[\rhola,\infty)$, the inclusions
\begin{equation}\label{pt9}
S_\rho(\lambda)\subset\{\bx\}+\ell\Vert\lambda-\bl\Vert\B\;\mbox{ and }\;\emp\ne S_\rho(\lm)\subset\inter\B_{\hat\gamma}(\bx),
\end{equation}
which mean that the mapping $S_\rho$ enjoys the strengthened robust isolated upper Lipschitzian property at $(\ox,\olm)$ uniformly in $\rho$ on the interval $[\rhola,\infty)$.
\end{theorem}\vspace*{-0.15in}
\begin{proof}
Since $\Phi$ is twice differentiable at $\bx$, there are constants $\hat\gamma\in(0,\gamma_{\olm}]$ and $\kappa>0$ with
\begin{equation}\label{eq49}
\Vert\Phi(x)-\Phi(\bx)\Vert\le\kappa\Vert x-\bx\Vert\;\textrm{ for all }\;x\in\B_{\hat\gamma}(\bx).
\end{equation}
Employing the second-order growth condition \eqref{eq16} tells us that $S_\rho(\bl)=\{\bx\}$ for all $\rho\ge\rhola$.
Define now the the positive constant
\begin{equation}\label{eq21}
\ell:=\dfrac{\kappa}{\ellla}+\sqrt{\frac{\kappa^2}{\ell_\bl^2}+\frac{1}{\ellla\rhola}},
\end{equation}
select a positive number $\epsilon<\ell^{-1}\hat\gamma$, and then pick any $\lambda\in\B_\epsilon(\olm)$ and $\rho\ge\rhola$. Observe further that for all such $\lm$ and $\rho$ we have $S_\rho(\lambda)\ne\emp$, since the optimization problem in \eqref{eq177} admits an optimal solution by the classical Weierstrass theorem. Fix any $u\in S_\rho(\lm)$ and recall from Proposition~\ref{paug}(ii) that the function $\lambda\mapsto\L(u,\lambda,\rho)$ is concave. This together with \eqref{eq15} yields
\begin{eqnarray}\label{eq18}
\L(u,\lambda,\rho)&\ge&\L(u,\bl,\rho)-\langle\nabla_\lambda\L(u,\lambda,\rho),\bl-\lambda\rangle\nonumber\\
&=&\L(u,\bl,\rho)-\rho^{-1}\left\langle\Pi_{\nQ}\big(\rho\Phi(u)+\lambda\big)-\lambda,\bl-\lambda\right\rangle\nonumber\\
&\ge& f(\bx)+\ellla\Vert u-\bx\Vert^2-{\rho^{-1}}\left\langle\Pi_{\nQ}\big(\rho\Phi(u)+\lambda\big)-\lambda,\bl-\lambda\right\rangle,
\end{eqnarray}
where we use \eqref{eq16} for the last inequality. It follows from the optimality of $u$ that
\begin{equation*}
\L(u,\lambda,\rho)\le\L(\bx,\lambda,\rho)=f(\bx)+\dfrac{\rho}{2}\dist^2\big(\Phi(\bx)+\rho^{-1}\lambda;\Q\big)-\dfrac{1}{2}\rho^{-1}\Vert\lambda
\Vert^2\le f(\bx),
\end{equation*}
which together with \eqref{eq18} brings us to the estimate
\begin{equation}\label{eq20}
\Vert u-\bx\Vert^2\le\dfrac{1}{\rho\ellla}\left\langle\Pi_{\nQ}\big(\rho\Phi(u)+\lambda\big)-\lambda,\bl-\lambda\right\rangle.
\end{equation}
Employing the projection properties (P2) and (P4) from Section~\ref{sect02}, we get
\begin{eqnarray*}
\big\Vert\Pi_{\nQ}\big(\rho\Phi(u)+\lambda\big)-\lambda\big\Vert
&=&\Vert\rho\Phi(u)-\Pi_\Q\big(\rho\Phi(u)+\lambda\big)\Vert\\
&=&\Vert\rho\left(\Phi(u)-\Phi(\bx)\right)+\Pi_\Q\big(\rho\Phi(\bx)+\bl\big)-\Pi_\Q\big(\rho\Phi(u)+\lambda\big)\Vert\\
&\le&\rho\Vert\Phi(u)-\Phi(\bx)\Vert+\rho\Vert\Phi(u)-\Phi(\bx)\Vert+\Vert\bl-\lambda\Vert\\
&\le&2\rho\kappa\Vert u-\bx\Vert+\Vert\bl-\lambda\Vert,
\end{eqnarray*}
where the last inequality comes from \eqref{eq49}. Using this and \eqref{eq20} tells us that
\begin{equation*}
\Vert u-\bx\Vert^2\le\dfrac{1}{\rho\ellla}\Big(2\rho\kappa\Vert u-\bx\Vert+\Vert\lambda-\bl\Vert\Big)\Vert\lambda-\bl\Vert,
\end{equation*}
which can be written in the equivalent form as
\begin{equation*}
\ellla\Vert u-\bx\Vert^2-2\kappa\Vert\lambda-\bl\Vert\cdot\Vert u-\bx\Vert-\dfrac{\Vert\lambda-\bl\Vert^2}{\rho}\le 0.
\end{equation*}
This in turn gives us the estimate
\begin{equation*}
\Vert u-\bx\Vert\le\left(\dfrac{\kappa}{\ellla}+\sqrt{\frac{\kappa^2}{\ell_\bl^2}+\frac{1}{\ellla\rho}}\right)\Vert\lambda-\bl\Vert\le\ell\Vert\lambda-\bl\Vert\le\ell\epsilon<\hat\gamma,
\end{equation*}
which simultaneously verifies both inclusions in \eqref{pt9} and thus completes the proof.
\end{proof}\vspace*{-0.05in}

It follows from Theorem~\ref{solv} that, at each iteration $k$, the condition $\lk\in\B_\ve(\bl)$ on the current multiplier in Algorithm~\ref{Alg1} allows us to find an exact local solution to the optimization problem \eqref{ag1} such that $\Vert u^k-\bx\Vert\le\ell\Vert\lk-\bl\Vert$. Then the Lipschitz continuity of $\nabla_x\L(\cdot,\lk,\rok)$ around $u^k$ ensures that for any $\ek\ge 0$ we can get an $\ek$-solution $\xkk$ satisfying both the approximate stationary condition \eqref{xkk} and the same estimate
\begin{equation}\label{isocal}
\Vert\xkk-\bx\Vert\le\ell\Vert\lk-\bl\Vert
\end{equation}
as the exact solution $u^k$ to the optimization problem \eqref{ag1} under consideration.\vspace*{0.05in}

Now we are ready to proceed with local convergence analysis of Algorithm~\ref{Alg1}, which mainly exploits the two major ingredients and the corresponding results developed above: {\bf(1)} {\em SOSC} \eqref{SOSC} at $(\ox,\olm)$ and the associated {\em second-order growth} of the augmented Lagrangian, and {\bf(2)} the {\em calmness} of the  multiplier mapping. In addition, we assume that the set of Lagrange multipliers is a {\em singleton}. The main reason for imposing this restriction is that the convergent analysis of the general case, namely when the set of Lagrange multipliers is not a singleton,  is conducted by using an iterative framework proposed by Fischer in \cite[Theorem~1]{Fis02}. However, the latter result demands an error bound estimate the for consecutive terms of the ALM. Deriving such an estimate for SOCPs is our ongoing research project.\vspace*{0.05in}

The following theorem establishes the linear convergence of Algorithm~\ref{Alg1} in both exact and inexact frameworks of the ALM with an arbitrarily chosen tolerance in \eqref{xkk} in the form $\ek=o\big(\sigma(\xk,\lk)\big)$, where $\sigma(x,\lm)$ is the error bound from \eqref{res}.\vspace*{-0.05in}

\begin{theorem}[\bf primal-dual convergence of ALM]\label{locpd} Let $(\bx,\bl)$ be a solution to the KKT system \eqref{KKT}, let SOSC \eqref{SOSC} hold at $(\ox,\olm)$, and let
the multiplier mapping $M_{\ox}$ from \eqref{mulmap} be calm at $((0,0),\olm)$ and $\Lm(\ox)=\{\olm\}$. Then there exist positive numbers $\bar\gamma$ and $\bar\rho$ ensuring the following: for any starting point $(x^0,\lambda^0)\in\B_{\bar\gamma}(\bx,\bl)$ and any $\rok\ge\bar\rho$, Algorithm~{\rm\ref{Alg1}} generates a sequence of iterates $(\xk,\lk)$ with a tolerance in \eqref{xkk} arbitrary chosen as $\ek=o(\sigma(\xk,\lk))$ such that $(\xk,\lk)$ converges to $(\ox,\bl)$ as $k\to\infty$, and the rate of this convergence is linear.
\end{theorem}\vspace*{-0.2in}
\begin{proof} Let $\rhola,\gamma_\bl$, $\ellla$ be the positive constants taken from Theorem~{\rm\ref{growth}(iii)}, and let $\kappa_i$ and $\gamma_i$ for $i=1,2,3$ be positive constants taken from the Lipschitzian estimates \eqref{eb}, \eqref{eb1}, and \eqref{pt1}, respectively. Picking the positive constants $\kappa$ and $\hat\gamma$ from \eqref{eq49}, $\ell$ from \eqref{eq21}, and $\epsilon$ from Theorem~\ref{solv}, define the positive numbers
\begin{equation}\label{cons1}
\hat\gamma_1:=\min\big\{\gamma_1,\hat\gamma\big\},\quad\gamma_{1,2}:=\max\big\{\gamma_1,\gamma_2\big\},\quad\gamma:=\min\Big\{\gamma_3,\dfrac{\hat\gamma_1}{2\kappa_3}\Big\},
\end{equation}
\begin{equation}\label{cons2}
\bar\rho:=\max\big\{\rhola,2\kappa_1,8\kappa_1^2\kappa_2\big\},\;\;\mbox{and}\;\;\bar\gamma:=\min\Big\{\hat\gamma_1,\gamma_2,\epsilon,\dfrac{\bar\rho\gamma}{2\sqrt{10}}, \dfrac{\hat\gamma_1}{2\ell},\dfrac{\gamma}{2\ell(\kappa+1)}\Big\}.
\end{equation}
Assume also without loss of generality that
\begin{equation}\label{ep}
o\big(\sigma(x,\lm)\big)\le \min\Big\{\dfrac{1}{\kappa_2\bar\rho},\frac{1}{8\kappa_1\kappa_2}\Big\}\sigma(x,\lm)\;\;\mbox{whenever}\;\; (x,\lm)\in \B_{\gamma_{1,2}}(\ox,\olm)
\end{equation}
and then show that for any starting point $(x^0,\lambda^0)\in\B_{\bar\gamma}(\bx,\bl)$ there exists a sequence $\{(\xk,\lk)\}$ generated by Algorithm~\ref{Alg1} with any $\rok\ge\bar\rho$ such that
\begin{equation}\label{pt10}
(\xk,\lk)\in\B_{\bar\gamma}(\bx,\bl)\;\mbox{for all }\;k\in\N\cup\{0\}.
\end{equation}
Arguing by induction, observe that \eqref{pt10} obviously holds for $k=0$ and suppose that \eqref{pt10} is satisfied for some $k\in\N$ with $\rok\ge\bar\rho$. We are going to verify that \eqref{pt10} fulfills for $k+1$. To furnish this, deduce first from \eqref{cons2} that $\|\lm^k-\olm\|\le\epsilon$. This together with the remark after the proof of Theorem~\ref{solv} ensures the existence of an approximate solution $\xkk$ with
\begin{equation*}
\Vert\nabla_x\L(\xkk,\lk,\rok)\Vert\le\ek=o\big(\sigma(\xk,\lk)\big),
\end{equation*}
where $\ve_k\ge 0$ can be chosen arbitrary in this form. It follows from \eqref{isocal} that the obtained $\ek$-solution satisfies the estimates
\begin{equation}\label{pt11}
\Vert\xkk-\bx\Vert\le\ell\|\lm^k-\olm\|\le\ell\bar\gamma\le\dfrac{\hat\gamma_1}{2},
\end{equation}
where the last inequality comes from \eqref{cons2}. We proceed now to establish a similar estimate for the dual iterate $\lm^{k+1}$. Using \eqref{lkk} and the projection property (P4) yields
$\lkk\in N_\Q(\Phi(\xkk)+\rho_k^{-1}(\lk-\lkk))$ and hence $\lkk\in M(\xkk,v^{k+1},w^{k+1})$ with $w^{k+1}:=\dfrac{\lk-\lkk}{\rok}$ and
\begin{equation}\label{pt6}
v^{k+1}:=\nabla_x L(\xkk,\lkk)=\nabla_x\L(\xkk,\lk,\rok)=o\big(\sigma(\xk,\lk)\big).
\end{equation}
The inclusion $(\xk,\lk)\in\B_{\gamma_2}(\bx,\bl)$ allows us to deduce from \eqref{eb1}, \eqref{ep}, and \eqref{pt6} that
\begin{equation}\label{pt4}
\Vert v^{k+1}\Vert\le\dfrac{\sigma(\xk,\lk)}{\kappa_2\bar\rho}\le\dfrac{\Vert\xk-\bx\Vert+\Vert\lk-\bl\Vert}{\bar\rho}.
\end{equation}
Employing again the updating scheme \eqref{lkk}, we arrive at the relationships
\begin{eqnarray*}\label{pt5}
\Vert w^{k+1}\Vert&=&\Vert\rho_k^{-1}\lk-\Pi_\nQ\big(\Phi(\xkk)+\rho_k^{-1}\lk\big)\Vert\nonumber\\
&\le&\rho_k^{-1}\Vert\lk-\bl\Vert+\big\Vert\Pi_\nQ\big(\Phi(\xkk)+\rho_k^{-1}\lk\big)-\Pi_\nQ\big(\Phi(\bx)+\rho_k^{-1}\bl\big)\big\Vert\nonumber\\
&\le&2\rho_k^{-1}\Vert\lk-\bl\Vert+\Vert\Phi(\xkk)-\Phi(\bx)\Vert\nonumber\\
&\le&2\bar\rho^{-1}\Vert\lk-\bl\Vert+\kappa\Vert\xkk-\bx\Vert
\end{eqnarray*}
with the last estimate coming from \eqref{eq49} and $\xkk\in\B_{\hat\gamma}(\bx)$. Thus \eqref{pt11} and \eqref{pt4} bring us to
\begin{eqnarray*}
\Vert\xkk-\bx\Vert+\Vert v^{k+1}\Vert+\Vert w^{k+1}\Vert&\le&(\kappa+1)\Vert\xkk-\bx\Vert+\bar\rho^{-1}\Vert\xk-\bx\Vert+3\bar\rho^{-1}\Vert\lk-\bl\Vert\\
&\le&\ell(\kappa+1)\Vert\lk-\bl\Vert+\sqrt{10}\bar\rho^{-1}\Vert(\xk,\lk)-(\bx,\bl)\Vert\le\gamma,
\end{eqnarray*}
where the last inequality employs the induction assumption \eqref{pt10} together with \eqref{cons2}. This along with $\gamma\le\gamma_3$ due to \eqref{cons1} ensures that $(\xkk,v^{k+1},w^{k+1})\in\B_{\gamma_3}(\bx,0,0)$. Hence we deduce from the upper Lipschitzian property in \eqref{pt1} and the definition of $\gamma$ in \eqref{cons1} that
\begin{equation*}\label{pt12}
\Vert\lkk-\bl\Vert\le\kappa_3\big(\Vert\xkk-\bx\Vert+\Vert v^{k+1}\Vert+\Vert w^{k+1}\Vert\big)\le\kappa_3\gamma\le\dfrac{\hat\gamma_1}{2}
\end{equation*}
verifying therefore the promised estimate for the dual iterate $\lm^{k+1}$. This together with \eqref{pt11} shows that $(\xkk,\lkk)\in\B_{\hat\gamma_1}(\bx,\bl)$. Using the latter, the imposed SOSC \eqref{SOSC}, and the calmness of the multiplier mappings $M_{\ox}$ from \eqref{mulmap}, we conclude from Theorem~\ref{error} that
\begin{equation*}
\Vert\xkk-\bx\Vert+\|\lkk-\olm\|\le\kappa_1\sigkk\;\mbox{ with}
\end{equation*}
\begin{equation}\label{eq7}
\sigkk=\Vert\nabla_x L(\xkk,\lkk)\Vert+\Vert\Phi(\xkk)-\Pi_\Q\big(\Phi(\xkk)+\lkk\big)\Vert.
\end{equation}
Define further the projection vector
\begin{equation*}
\skk:=\Pi_\Q(\Phi(\xkk)+\rho_k^{-1}\lk)
\end{equation*}
and deduce from the updating scheme \eqref{lkk} that
\begin{equation}\label{eq3}
\Phi(\xkk)-\skk=\dfrac{\lkk-\lk}{\rok}.
\end{equation}
Employing the projection properties (P1) and (P2) results in $\langle\skk,\lkk\rangle=0$ due to
\begin{equation*}\label{eq1}
{\rho_k^{-1}}{\lkk}=\Pi_\nQ\big(\Phi(\xkk)+{\rho_k^{-1}}{\lk}\big)=\Phi(\xkk)+{\rho_k^{-1}}{\lk}-\skk,
\end{equation*}
which together with $\skk\in\Q$ yields $\lkk\in N_\Q(\skk)$. Hence $\skk=\Pi_\Q(\skk+\lkk)$ by property (P4). Since the mapping $y\mapsto y-\Pi_\Q(y+\lkk)=\Pi_\nQ(y+\lkk)-\lkk$ is clearly nonexpansive, we arrive at the relationships
\begin{eqnarray*}
&&\big\Vert\Phi(\xkk)-\Pi_\Q\big(\Phi(\xkk)+\lkk\big)\big\Vert\nonumber\\
&=&\big\Vert\Phi(\xkk)-\Pi_\Q\big(\Phi(\xkk)+\lkk\big)\big\Vert-\big\Vert\skk-\Pi_\Q(\skk+\lkk)\big\Vert\nonumber\\
&\le&\big\Vert\Phi(\xkk)-\Pi_\Q\big(\Phi(\xkk)+\lkk\big)-\big(\skk-\Pi_\Q(\skk+\lkk)\big)\big\Vert\nonumber\\
&\le&\Vert\Phi(\xkk)-\skk\Vert\\
&\le&\rho_k^{-1}\Vert \lkk-\lk\Vert\quad \quad(\mbox{by}\;\;\eqref{eq3})\\
&\le&\rho_k^{-1}\big(\Vert\lkk-\bl\Vert+\Vert\lk-\bl\Vert\big)\\
&\le&\kappa_1\rho_k^{-1}(\sigkk+\sigk).
\end{eqnarray*}
Using this together with \eqref{pt6} and \eqref{eq7} leads us to the estimates
\begin{equation*}
\sigkk\le\ek+\big\Vert\Phi(\xkk)-\Pi_\Q\big(\Phi(\xkk)+\lkk\big)\big\Vert\le\ek+\dfrac{\kappa_1}{\rok}\Big(\sigkk+\sigk\Big),
\end{equation*}
which can be equivalently rewritten as
\begin{equation*}\label{eq8}
\Big(1-\dfrac{\kappa_1}{\rok}\Big)\sigkk\le\ek+\dfrac{\kappa_1}{\rok}\sigk.
\end{equation*}
Since $\rok\ge\bar\rho$, by \eqref{cons2}, we get  $1-\dfrac{\kappa_1}{\rok}>\dfrac{1}{2}$, which  ensures  that
\begin{equation*}
\sigkk\le 2\sigk\left(\dfrac{\ek}{\sigk}+\dfrac{\kappa_1}{\rok}\right).
\end{equation*}
Applying finally the error bounds \eqref{eb} and \eqref{eb1} and then appealing to \eqref{cons2} and \eqref{ep} yields
\begin{eqnarray}\label{eq11}
\Vert\xkk-\bx\Vert+\Vert\lkk-\bl\Vert&\le&\kappa_1\sigkk\le 2\kappa_1\left(\dfrac{\ek}{\sigk}+\dfrac{\kappa_1}{\rok}\right)\sigk\nonumber\\
&\le&2\kappa_1\kappa_2\left(\dfrac{\ek}{\sigk}+\dfrac{\kappa_1}{\rok}\right)\big(\Vert\xk-\bx\Vert+\Vert\lk-\bl\Vert\big)\nonumber\\
&\le&\frac{1}{2}\big(\Vert\xk-\bx\Vert+\Vert\lk-\bl\Vert\big),
\end{eqnarray}
which together with the induction assumption \eqref{pt10} brings us to
\begin{equation*}
(\xkk,\lkk)\in\B_{\bar\gamma}(\bx,\bl).
\end{equation*}
This finishes our induction argument to justify \eqref{pt10} for all $k\in\N$. Observe that the latter inclusion along with \eqref{cons2} implies that $\|\lkk-\olm\|\le\epsilon$ while allowing us to use Theorem~\ref{solv} to construct the next primal iterate $x^{k+2}$. Since \eqref{eq11} holds for all $k\in \N$, we clearly get that $(x^k,\lm^k)\to(\ox,\olm)$ as $k\to\infty$.
Furthermore, the obtained estimate tells us that rate of convergence of $(x^k,\lm^k)$ to $(\ox,\olm)$ is linear, which therefore completes the proof of the theorem.
\end{proof}\vspace*{-0.05in}

To conclude the paper, let us compare the convergence analysis of the ALM given in Theorem~\ref{locpd} with the one provided recently by Kanzow and Steck \cite{KaS17,KaS19} for the class of
${\cal C}^2$-cone reducible conic programs that includes SOCPs. These publications were devoted to convergence analysis of  a modified version of the ALM, called the {\em safeguarded augmented
Lagrangian methods} therein, in which the $\lm^k$  in the formation of subproblems \eqref{ag1} is replaced with a certain vector $w^k$ chosen from a bounded set or obtained from the projection of
$\lm^k$ onto a bounded set. The main motivation for this modification comes from the fact that  the dual sequence, constructed by the standard ALM, may be unbounded in general. As has been
extensively documented in \cite{abm1,abm3,abm2,bm14}, this modification has  a remarkable  global convergence theory and was successfully implemented in the ALGENCAN software; see
\cite[Chapter~10]{bm14} for a through discussion about the different aspects of this implementation. Also, note that while the safeguarded ALM in \cite{KaS17,KaS19} uses a particular updating
scheme for the penalty parameter $\rho_k$, our approach reveals that there is no need to confine the convergence analysis of the standard ALM to a particular updating scheme for $\rho_k$ as long
as we keep it sufficiently large. Also, as mentioned in Section~\ref{intro}, the solvability of the standard ALM subproblems  \eqref{ag1} was not addressed in \cite{KaS17,KaS19}. Let us finally
emphasize that the progress achieved in this paper is largely based on the application and development of powerful tools of second-order variational analysis and generalized
differentiation.
\end{document}